\documentclass[12pt,letterpaper]{article}

\UseRawInputEncoding
\usepackage{bbm}
\usepackage{chngcntr}
\usepackage{multirow}

\usepackage[colorlinks=true,urlcolor=blue,citecolor=blue,linkcolor=blue,bookmarks=true,bookmarksopen=false]{hyperref}

\usepackage{algorithmic}
\usepackage[linesnumbered,ruled,vlined]{algorithm2e}

\usepackage[left=1in,top=1in,right=1in,bottom=1in]{geometry}

\newcommand{\exclude}[1]{}

\renewcommand{\Omega}{\varOmega}

\newcounter{commentcounter}
\long\def\symbolfootnote[#1]#2{\begingroup%
\def\thefootnote{\fnsymbol{footnote}}\footnote[#1]{#2}\endgroup}

\usepackage{sibarticle}

\usepackage{amsmath}
\usepackage{kbordermatrix}
\usepackage{mathtools}

\usepackage{lineno}

\usepackage{tikz}
\usetikzlibrary{arrows,shapes,snakes,automata,backgrounds,petri}
\usetikzlibrary{graphs}

\newcommand{\ignore}[1]{}

\usepackage{color}
\usepackage{bm}
\usepackage{multirow}
\usepackage{booktabs}

\usepackage{rotating}
\usepackage[symbol]{footmisc}
\usepackage{graphicx}
\usepackage{subfigure}
\usepackage{epstopdf}

\allowdisplaybreaks

\title{On Exact Dual and Mixed-Integer Reformulations for Finite Horizon Markov Decision Processes with Submodular Revenue Maximization}
\ShortTitle{MIP for MDPs with Submodular Revenue Maximization} \ShortAuthors{Authors}
\NumberOfAuthors{1} \FirstAuthor{Jia-Yu Chen, Hsin-Yuan Huang, Tzu-Yang Teng, Yun-Jen Tsai, Hao-Hsiang Wu\footnotemark[1]}
\FirstAuthorAddress{ Department of Management Science, National Yang Ming Chiao Tung University, Hsinchu, Taiwan\\ 
	{\tt \{cjy06.mg13, w998688.mg11, asd304go.mg10, yzc.mg14,  hhwu2\}@nycu.edu.tw	}}

\SecondAuthor{Authors}
\SecondAuthorAddress{
	TBD\\ {\tt
		TBD}
}
\ShortAuthors{Chen, Huang, Teng, Tsai, and Wu}
\exclude{

}

\renewcommand{\baselinestretch}{1.5}

\keywords{Mixed-Integer Programming; Duality;  Submodularity; Markov Decision Process;}

\begin{document}

\maketitle \centerline{\today}


\footnotetext[1]{Corresponding author (hhwu2@nycu.edu.tw)} 
\footnotetext[0]{}
\begin{abstract}
We study exact mixed-integer reformulations of a finite horizon Markov decision process (MDP) with transition revenues defined by local submodular revenue maximization problems. We start with the dual linear programming formulation of the MDP, where local
optimal revenues appear in the objective function and the constraints
describe the state probabilities from the initial state. We combine the dual formulation with the binary variables of the
local problems and obtain an exact mixed-integer nonlinear
formulation with bilinear terms. We use the submodular inequalities and an exact McCormick linearization to obtain an equivalent mixed-integer linear
formulation. We also use the dual structure to conduct a sensitivity analysis of
the MDP optimal value with respect to changes in the local revenue
coefficients. Based on local incumbent revenues and upper bounds, this sensitivity
analysis evaluates the effect of local optimality gaps on the global
MDP optimality gap, without solving every local revenue maximization
problem to optimality.

\end{abstract}

\renewcommand{\baselinestretch}{1.5}

\section{Introduction}\label{sec:MDP_LP}

We consider exact mixed-integer reformulations for a finite horizon Markov decision process (MDP) with embedded local submodular revenue maximization problems. The optimal values of these embedded problems determine the revenues that are associated with the possible transitions. Dynamic decision models use designed reward and cost functions to represent the consequences of decisions in various applications \cite{AlagozMS2004,ElwanyOR2011,HoppWuIIE1990,KalymonMS1972,LiaoEJOR2006,ShechterOR2008,VeatchWeinOR1994,ZuckermanNRL1986}. These representative applications include, but are not limited to, maintenance, production, inventory control, and several medical treatment decisions. We consider a different revenue structure, where a separate local submodular optimization problem determines the revenue associated with each possible transition. Our model is inspired by the stochastic depletion framework of Chan and Farias \cite{ChanMOR2009}. In what follows, we adopt the notation standards from Puterman \cite{Puterman1994} to introduce the MDP model for consistency with the classic MDP literature, where decisions occur at discrete epochs, the horizon is finite, and parameters may depend on time. We also use the terminology and  examples from \cite{ChanMOR2009} to introduce the MDP components, without adding further restrictions to the  MDP model. 

Formally, we define $T\in\mathbb{N}$ as the terminal epoch and simply define the set of decision epochs as $\mathbb{T}=\{0,1,\ldots,T-1\}$. Let $\mathbb{S}$ be a finite state space, and let $\{S_t\}_{t=0}^T$ denote the stochastic state process, where each $S_t$ takes values in $\mathbb{S}$. In the stochastic depletion framework of \cite{ChanMOR2009}, each state records the available items, the current epoch, and information that describes the depletion process. We use the general state $s\in\mathbb{S}$ and apply no specific depletion structure on its features. Besides, let $\mathbb{A}$ be an action space and let $\{A_t\}_{t=0}^{T-1}$ denote the stochastic action process, where each $A_t$ takes values in $\mathbb{A}$. For each decision epoch $t\in\mathbb{T}$ and state $s\in\mathbb{S}$, let $\mathcal{A}_t(s)\subseteq\mathbb{A}$ denote the finite and nonempty set of actions that are available in state $s$ at epoch $t$. An activity in the stochastic depletion framework plays the role of an action $a\in\mathcal{A}_t(s)$, where the selected activity determines the random depletion of the available items. In the dynamic product line design application of \cite{ChanMOR2009}, an activity corresponds to a feasible product line under a cardinality constraint. We next specify the transition and revenue structures of the MDP model. Related studies use reward and transition structures to analyze value functions and optimal decisions
\cite{LightMOR2021,SmithMcCardleOR2002}. Given a decision epoch $t\in\mathbb{T}$, a state $s\in\mathbb{S}$, and an action $a\in\mathcal{A}_t(s)$, the transition probability from state $s$ to state $s'\in\mathbb{S}$ is defined as $P_t(s'\mid s,a)=\Pr(S_{t+1}=s'\mid S_t=s,A_t=a)$ and $\sum_{s'\in\mathbb{S}}P_t(s'\mid s,a)=1$. In the stochastic depletion framework of \cite{ChanMOR2009}, the random depletion that follows the selected activity induces the transition probabilities. We instead treat $\{P_t\}_{t\in\mathbb{T}}$ as part of the MDP data and impose no depletion structure on these probabilities. We define $r_t(s,a,s')\in\mathbb{R}$ as the immediate revenue that is obtained when action $a$ is selected in state $s$ at epoch $t$ and the next state is $s'$. In the submodular stochastic depletion model of \cite{ChanMOR2009}, the submodularity applies to the cumulative value of depleted items. The reward obtained during each period equals the corresponding marginal increase in this cumulative value. Related MDP studies use submodularity conditions on immediate reward or cost functions defined over ordered state and action spaces \cite{AltmanQS1995,KilicEJOR2026, LuOR1984}. Another study defines the stage replenishment cost as a submodular function of an action represented by a subset \cite{ZhangMSOM2018}. In these models, the MDP action enters a prescribed reward or cost function directly.

Combinatorial optimization arises in the construction of the myopic policy in \cite{ChanMOR2009}. At each decision epoch, this policy solves an optimization problem to select an activity that maximizes the expected immediate reward. We use these observations in a different construction. We study a separate local problem of constrained submodular maximization for each possible transition $(t,s,a,s')$. The local solution does not represent the MDP action. Also, its optimal value defines the transition revenue $r_t(s,a,s')$. Therefore, the submodularity is for the local decision variables $\mathbf{x}^{(t,s,a,s')}$ rather than the action $a\in\mathcal{A}_t(s)$. Motivated by this use of submodularity, we assume that $r_t(s,a,s')$ is the optimal value of a local submodular revenue maximization problem associated with $(t,s,a,s')$. That is, for each $t\in\mathbb{T}$, $s,s'\in\mathbb{S}$, and $a\in\mathcal{A}_t(s)$, let $n^{(t,s,a,s')}\in \mathbb{N}$ denote the number of local decision variables. Here, we define the ground set as $\mathcal{I}^{(t,s,a,s')}=\{1,\ldots,n^{(t,s,a,s')}\}$. We use $\mathbf{x}^{(t,s,a,s')}\in\mathbb{B}^{n^{(t,s,a,s')}}$ to represent the binary decision vector of a subset of $\mathcal{I}^{(t,s,a,s')}$. Let $f^{(t,s,a,s')}:2^{\mathcal{I}^{(t,s,a,s')}}\rightarrow\mathbb{R}$ denote the corresponding submodular revenue function, where $f^{(t,s,a,s')}(U)+f^{(t,s,a,s')}(W)
\geq
f^{(t,s,a,s')}(U\cup W)+f^{(t,s,a,s')}(U\cap W)
$ for all $U,W\subseteq\mathcal{I}^{(t,s,a,s')}$. Let $\mathcal{X}^{(t,s,a,s')}\subseteq\mathbb{R}^{n^{(t,s,a,s')}}$ denote the feasible region defined by a set of linear constraints of the local problem. We suppose that $\mathcal{X}^{(t,s,a,s')}\cap\mathbb{B}^{n^{(t,s,a,s')}}$ is nonempty  and finite. With a slight abuse of notation, we write $f^{(t,s,a,s')}(\mathbf{x}^{(t,s,a,s')})$ for the value of the subset that is represented by $\mathbf{x}^{(t,s,a,s')}$. For any
$\bar{\mathbf{x}}^{(t,s,a,s')}
\in \mathcal{X}^{(t,s,a,s')}\cap \mathbb{B}^{n^{(t,s,a,s')}}$,
let
$\overline{X}^{(t,s,a,s')}
=
\{i\in\mathcal{I}^{(t,s,a,s')}:
\bar{x}_i^{(t,s,a,s')}=1\}$
denote its support. For notational convenience, we use
$\bar{\mathbf{x}}^{(t,s,a,s')}$ and
$\overline{X}^{(t,s,a,s')}$ interchangeably. We also use
$f^{(t,s,a,s')}(\bar{\mathbf{x}}^{(t,s,a,s')})$ and
$f^{(t,s,a,s')}(\overline{X}^{(t,s,a,s')})$ interchangeably for the
same function value. We call $\overline{X}^{(t,s,a,s')}$ feasible for
the local problem if
$\bar{\mathbf{x}}^{(t,s,a,s')}
\in
\mathcal{X}^{(t,s,a,s')}
\cap
\mathbb{B}^{n^{(t,s,a,s')}}$. For each $t\in\mathbb{T}$, $s,s'\in\mathbb{S}$, and $a\in\mathcal{A}_t(s)$, the transition revenue on the possible transition $(t,s,a,s')$ is defined as a local submodular revenue maximization problem
\begin{equation}\label{MDP_transition_revenue}
r_t(s,a,s')
=
\max
\left\{
f^{(t,s,a,s')}(\mathbf{x}^{(t,s,a,s')}):
\mathbf{x}^{(t,s,a,s')}
\in
\mathcal{X}^{(t,s,a,s')}
\cap
\mathbb{B}^{n^{(t,s,a,s')}}
\right\}.
\end{equation}
Since the local feasible set $\mathcal{X}^{(t,s,a,s')}
\cap
\mathbb{B}^{n^{(t,s,a,s')}}$ is finite and nonempty, the corresponding optimal solution of the local problem \eqref{MDP_transition_revenue} exists and is finite.
The submodular revenue function $f^{(t,s,a,s')}$ and the feasible region $\mathcal{X}^{(t,s,a,s')}$ may vary across decision epochs, current states, actions, and next states. The local submodular function maximization problem \eqref{MDP_transition_revenue} is NP-hard, even for monotone functions subject to a cardinality constraint \cite{Feige1998, Feige2011, NWF78}. An exact integer programming formulation for the submodular maximization problem
\eqref{MDP_transition_revenue} is introduced by Nemhauser and  Wolsey \cite{NW81}.
This formulation is based on submodular inequalities and has provided
a foundation for subsequent polyhedral studies \cite{Ahmed2011,Coniglio2022,Shi2022,first2016,Yu2017}. Here, we refer
readers to the tutorial in \cite{KYTut2023} for a general review of this literature. Also, related polyhedral methods based on submodularity appear in chance-constrained optimization  \cite{AtamturkDO2009,AtamturkBhardwaj2015, Fatma2022, Shen2022, Wu2017} and conic quadratic
optimization \cite{Andres2018,Alper2020,Alper2022,kilincc2020conic}.
Similar developments have been established for $k$-submodular
optimization \cite{QYuk2021,QYu2021} and submodular minimization \cite{YuAhmed2017,QYu2023,YuQ2022}. More recent work extends this polyhedral perspective to functions satisfying $L^\natural$ convexity \cite{YuKucukyavuz2025} and robust mixed-integer programs \cite{WuIJOO2026}. This line of research motivates our study of the connection
between the present MDP model and the polyhedral structure defined by
submodular inequalities.

We provide a product line interpretation for another view of our model. In a market or manufacturing environment, the MDP action $a\in\mathcal{A}_t(s)$ represents a strategic market or manufacturing capacity decision. The selected action affects the subsequent evolution of the market or manufacturing state through a transition $P_t(\cdot\mid s,a)$. After the transition from state $s$ to state $s'$, the  associated firm selects a feasible product portfolio $\mathbf{x}^{(t,s,a,s')}$ under the corresponding market and operational conditions. However, adding a new product does not necessarily generate significant additional profits. In this product line interpretation, the submodularity may be interpreted as diminishing marginal contribution to revenue from each additional product in the portfolio. That is, feasibility requires $\mathbf{x}^{(t,s,a,s')}\in\mathcal{X}^{(t,s,a,s')}\cap\mathbb{B}^{n^{(t,s,a,s')}}$. The index $(t,s,a,s')$ identifies these conditions and the associated local revenue submodular maximization problem. Here, we note that the optimal value is shown as $r_t(s,a,s')$, where the local decision variables affect only the revenue that is earned on the realized transition. These variables do not enter $P_t(\cdot\mid s,a)$ and have no effect on subsequent state evolution. Therefore, their optimal value can be represented by $r_t(s,a,s')$ without enlarging the MDP action or state space. Also, let $r_T(s)\in\mathbb{R}$ be the terminal revenue that is obtained when the process terminates in state $s$ at epoch $T$. The authors of \cite{ChanMOR2009} restrict the depletion rewards to nonnegative values and set the terminal reward to zero. We allow the stage and terminal revenues to take constant values in $\mathbb{R}$. We assume that the stage and terminal revenue functions are bounded, where a constant $R_{\max}$ has $r_t(s,a,s')\leq R_{\max}$ and $r_T(s)\leq R_{\max}$ for all $t\in\mathbb{T}$, $s,s'\in\mathbb{S}$, and $a\in\mathcal{A}_t(s)$. Finally, the finite horizon MDP with the embedded local submodular revenue maximization problems can be represented as a tuple
\begin{equation}\label{MDP_model}
\mathcal{M}=\left(\mathbb{S},\mathbb{A},T,
\{\mathcal{A}_t(s)\}_{t\in\mathbb{T},\,s\in\mathbb{S}},
\{P_t\}_{t\in\mathbb{T}},
\{r_t\}_{t\in\mathbb{T}},
r_T \right).
\end{equation}
Here, the MDP model \eqref{MDP_model} contains two optimization levels, where the local level determines the transition revenue coefficients and the MDP level determines a policy. In the next section, we study this problem through the dual linear programming formulation of the MDP model \eqref{MDP_model}. Based on this dual structure, we develop exact mixed-integer reformulations that integrate both levels within a single mathematical programming model. Note that, in the dual objective, each local revenue is weighted by a transition probability and its associated dual variable. In a subsequent section, we also use this relation to analyze the sensitivity of the optimal value of the MDP model \eqref{MDP_model} to changes in the local revenues. We use the incumbent revenues and upper bounds obtained for the local submodular
maximization problems \eqref{MDP_transition_revenue} to study bounds on the resulting optimality gap. This sensitivity analysis evaluates the effect of variations in the local revenue coefficients on the MDP optimal value.

\section{Exact Dual and Mixed-Integer Reformulations}

We begin with the classical LP formulation of finite horizon MDPs
\cite{DermanKlein1965,Manne1960,Puterman1994}. Its dual formulation provides the
structure used in our analysis. We first introduce the policy and
value functions needed to state this formulation. A randomized Markov policy is denoted by $\mathbf{\pi}=\{\pi_t\}_{t\in\mathbb{T}}$, where for each decision epoch $t\in\mathbb{T}$ and state $s\in\mathbb{S}$, the notation $\pi_t(\cdot\mid s)$ is a probability distribution on $\mathcal{A}_t(s)$. That is, $\pi_t(a\mid s)\geq 0$ for all $a\in\mathcal{A}_t(s)$ and $\sum_{a\in\mathcal{A}_t(s)}\pi_t(a\mid s)=1$. Unlike a stationary policy for an infinite horizon, the decision rule $\pi_t$ is allowed to depend explicitly on the decision epoch. Given an initial state $S_0=s$ and a policy $\mathbf{\pi}$, the finite horizon cumulative revenue is defined as
$G_0^{\mathbf{\pi}}
=
\sum_{t=0}^{T-1}r_t(S_t,A_t,S_{t+1})+r_T(S_T).$
In this finite horizon formulation, we consider the undiscounted case, which is equivalent to setting the discount factor to one. Here, we have $G_0^{\mathbf{\pi}}\leq(T+1)R_{\max}$. Thus, the state value function associated with $\mathbf{\pi}$ at every epoch $t\in\mathbb{T}$ is given by
\begin{equation*}\label{MDP_value}
V_t^{\mathbf{\pi}}(s)
=
\mathbb{E}_{\mathbf{\pi}}
\left[
\left.
\sum_{k=t}^{T-1}r_k(S_k,A_k,S_{k+1})+r_T(S_T)
\,\right|\,S_t=s
\right]
\end{equation*}
with the terminal condition $V_T^{\mathbf{\pi}}(s)=r_T(s)$ for all $s\in\mathbb{S}$. For each policy $\mathbf{\pi}$, the value functions of $V_t^{\mathbf{\pi}}(s)$ satisfy the Bellman expectation equations
\begin{equation*}\label{MDP_Bellman_expectation}
	V_t^{\mathbf{\pi}}(s)
	=\sum_{a\in\mathcal{A}_t(s)}\pi_t(a\mid s)
	\sum_{s'\in\mathbb{S}}P_t(s'\mid s,a)
	\left[r_t(s,a,s')+V_{t+1}^{\mathbf{\pi}}(s')\right]
	\quad \forall s\in\mathbb{S},\ t\in\mathbb{T}.
\end{equation*}
Since each set $\mathcal{A}_t(s)$ is finite and nonempty, the backward induction provides an optimal Markov policy $\mathbf{\pi}^*$. For every randomized Markov policy $\mathbf{\pi}$, this policy satisfies
$
V_t^{\mathbf{\pi}^*}(s)
\geq
V_t^{\mathbf{\pi}}(s)
$ for all $s\in\mathbb{S}$ and  $t\in\mathbb{T}$.
We define the optimal value function as $V_t^*(s)=V_t^{\mathbf{\pi}^*}(s)$ for all $s\in\mathbb{S}$ and $t\in\mathbb{T}$. The terminal condition is $V_T^*(s)=r_T(s)$ for all $s\in\mathbb{S}$. Therefore, the optimal value functions satisfy the Bellman optimality equations.
\begin{equation*}\label{MDP_Bellman_optimality}
V_t^*(s)
=
\max_{a\in\mathcal{A}_t(s)}
\left\{
\sum_{s'\in\mathbb{S}}P_t(s'\mid s,a)
\left[r_t(s,a,s')+V_{t+1}^*(s')\right]
\right\}
\quad \forall s\in\mathbb{S},\ t\in\mathbb{T}.
\end{equation*}

Here, the function $V_t^*(s)$ is the greatest expected revenue obtained for all feasible actions $a\in\mathcal{A}_t(s)$. This observation leads to an LP formulation of the finite horizon MDP optimization problem. Let $s_0\in\mathbb{S}$ be a fixed initial state and introduce the variables $v_t(s)\in\mathbb{R}$ for all $t\in\{0,1,\ldots,T\}$ and $s\in\mathbb{S}$. The finite horizon MDP optimization problem can then be formulated as follows. We note that the dual variables, $y_t(s,a)\geq 0$ for all $a\in\mathcal{A}_t(s),\
s\in\mathbb{S},\
t\in\mathbb{T}$ and $y_T(s)\in\mathbb{R}$ for $s\in\mathbb{S}$, are shown in brackets at each constraint so the reader can easily refer to them later for the corresponding dual formulation.
\begin{equation}\label{MDP_primal_LP}
\begin{aligned}
\min \quad
& v_0(s_0)
\\
\text{s.t.}\quad
& v_t(s)
\geq
\sum_{s'\in\mathbb{S}}
P_t(s'\mid s,a)
\left[
r_t(s,a,s')+v_{t+1}(s')
\right]
&&
\forall a\in\mathcal{A}_t(s),\
s\in\mathbb{S},\
t\in\mathbb{T},
&&
\bigl[y_t(s,a)\geq 0\bigr]
\\
& v_T(s)=r_T(s)
&&
\forall s\in\mathbb{S},
&&
\bigl[y_T(s)\in\mathbb{R}\bigr]
\\
& v_t(s)\in\mathbb{R}
&&
\forall s\in\mathbb{S},\
t\in\{0,1,\ldots,T\}.
&&
\end{aligned}
\end{equation}
The following argument briefly shows that the formulation \eqref{MDP_primal_LP} is an exact LP formulation for the MDP. Here, the constraints of the formulation \eqref{MDP_primal_LP} require $v_t(s)$ to dominate the expected revenue associated with every feasible action. By the backward induction, every feasible solution satisfies $v_t(s)\geq V_t^*(s)$ for all $s\in\mathbb{S}$ and $t\in\{0,1,\ldots,T\}$, including $v_0(s_0)\geq V_0^*(s_0)$. The optimal value functions $V_t^*(s)$ for all $s\in\mathbb{S},\ t\in\mathbb{T}$  satisfy the Bellman optimality equations. Thus, the optimal value functions are a feasible solution to the formulation \eqref{MDP_primal_LP}. It follows that the optimal objective value of the formulation \eqref{MDP_primal_LP} is $V_0^*(s_0)$.  

The duality of the formulation \eqref{MDP_primal_LP} provides the structure used in our subsequent analysis. The dual LP is given by
\begin{equation}\label{MDP_dual_LP}
\begin{aligned}
\max\quad
&\sum_{t\in\mathbb{T}}\sum_{s\in\mathbb{S}}\sum_{a\in\mathcal{A}_t(s)}
\sum_{s'\in\mathbb{S}}P_t(s'\mid s,a)r_t(s,a,s')y_t(s,a)
+\sum_{s\in\mathbb{S}}r_T(s)y_T(s)\\
\text{s.t.}\quad
&\sum_{a\in\mathcal{A}_0(s_0)}y_0(s_0,a)=1,\\
&\sum_{a\in\mathcal{A}_0(s)}y_0(s,a)=0
&&\forall s\in\mathbb{S}\setminus\{s_0\},\\
&\sum_{a\in\mathcal{A}_t(s)}y_t(s,a)
-\sum_{\tilde{s}\in\mathbb{S}}\sum_{a\in\mathcal{A}_{t-1}(\tilde{s})}
P_{t-1}(s\mid \tilde{s},a)y_{t-1}(\tilde{s},a)=0
&&\forall s\in\mathbb{S},\ t\in\mathbb{T}\setminus\{0\},\\
&y_T(s)-\sum_{\tilde{s}\in\mathbb{S}}\sum_{a\in\mathcal{A}_{T-1}(\tilde{s})}
P_{T-1}(s\mid \tilde{s},a)y_{T-1}(\tilde{s},a)=0
&&\forall s\in\mathbb{S},\\
&y_t(s,a)\in \mathbb R_+,
&&\forall a\in\mathcal{A}_t(s),\ s\in\mathbb{S},\ t\in\mathbb{T},\\
& y_T(s) \in \mathbb R,
&&\forall s\in\mathbb{S}.
\end{aligned}
\end{equation}
From the strong duality theorem (see, e.g., Chapter 4 of \cite{BT1997} and Chapter 7 of \cite{Schrijver98}), the primal and dual formulations \eqref{MDP_primal_LP} and \eqref{MDP_dual_LP} have the same optimal value $V_0^*(s_0)$. Note that mathematical programming formulations and duality can be used to analyze different structures in MDPs
\cite{BertsimasNinoMoraOR2000,LiShapiro2025,Ruszczynski2010} and to
derive bounds for stochastic dynamic programs based on dual
relaxations
\cite{AdelmanMersereauOR2008,BrownSmithSunOR2010}. In what follows,
motivated by these studies, we use the dual formulation 
\eqref{MDP_dual_LP} to analyze an MDP with embedded local
submodular revenue maximization problems \eqref{MDP_transition_revenue}. In our setting, the local
optimal revenues of the local problems \eqref{MDP_transition_revenue} appear as
coefficients in the dual objective function, and the feasible region
of the MDP dual variables remains unchanged. Based on this dual structure,
we obtain an exact mixed-integer bilinear formulation and an
equivalent mixed-integer linear formulation, and derive optimality gap bounds from the incumbent revenues and the upper bounds
obtained for the local problems 
\eqref{MDP_transition_revenue}.

Now, based on the formulation \eqref{MDP_dual_LP}, each transition revenue $r_t(s,a,s')$ appears as an objective coefficient. We next represent each coefficient through the submodular inequalities associated with the local problem \eqref{MDP_transition_revenue}. For each $S\subseteq\mathcal{I}^{(t,s,a,s')}$ and $i\in\mathcal{I}^{(t,s,a,s')}\setminus S$, we define the marginal contribution to revenue as $\rho_i^{(t,s,a,s')}(S)
=
f^{(t,s,a,s')}(S\cup\{i\})-f^{(t,s,a,s')}(S)$. Given $U\subseteq W\subseteq\mathcal{I}^{(t,s,a,s')}\setminus\{i\}$, the submodularity implies $\rho_i^{(t,s,a,s')}(U)
\geq
\rho_i^{(t,s,a,s')}(W)$.
Thus, the marginal revenue from an additional element cannot increase as the selected subset expands. We introduce the continuous variable $\theta^{(t,s,a,s')}\in\mathbb{R}$ to represent the approximation of the local revenue  $r_t(s,a,s')$. Given a  subset $\overline{X}^{(t,s,a,s')}\subseteq\mathcal{I}^{(t,s,a,s')}$, the corresponding submodular inequality of \cite{NW81} is
\begin{equation}\label{submodular_inequality}
\begin{aligned}
\theta^{(t,s,a,s')}
\leq
f^{(t,s,a,s')}(\overline{X}^{(t,s,a,s')})
&-
\sum_{i\in\overline{X}^{(t,s,a,s')}}
\rho_i^{(t,s,a,s')}(\mathcal{I}^{(t,s,a,s')}\setminus\{i\})
(1-x_i^{(t,s,a,s')})\\
&+
\sum_{i\in\mathcal{I}^{(t,s,a,s')}\setminus\overline{X}^{(t,s,a,s')}}
\rho_i^{(t,s,a,s')}(\overline{X}^{(t,s,a,s')})
x_i^{(t,s,a,s')}.
\end{aligned}
\end{equation}
The subset $\overline{X}^{(t,s,a,s')}$ indexes the inequality and need not satisfy the constraints of the local problem. For each $t\in\mathbb{T}$, $s,s'\in\mathbb{S}$, and $a\in\mathcal{A}_t(s)$, the submodular polyhedron of \cite{NW81} is defined as
\begin{equation*}\label{submodular_polyhedron}
\begin{aligned}
\mathcal{P}^{(t,s,a,s')}
=
\Bigg\{
(\mathbf{x}^{(t,s,a,s')},\theta^{(t,s,a,s')})
\in \mathbb B^{n^{(t,s,a,s')}}\times\mathbb{R}:\quad
&\theta^{(t,s,a,s')}
\leq
f^{(t,s,a,s')}(\overline{X}^{(t,s,a,s')})\\
&-
\sum_{i\in\overline{X}^{(t,s,a,s')}}
\rho_i^{(t,s,a,s')}(\mathcal{I}^{(t,s,a,s')}\setminus\{i\})
(1-x_i^{(t,s,a,s')})\\
&+
\sum_{i\in\mathcal{I}^{(t,s,a,s')}\setminus\overline{X}^{(t,s,a,s')}}
\rho_i^{(t,s,a,s')}(\overline{X}^{(t,s,a,s')})
x_i^{(t,s,a,s')}\\
&\forall\,
\overline{X}^{(t,s,a,s')}
\subseteq
\mathcal{I}^{(t,s,a,s')}
\Bigg\}.
\end{aligned}
\end{equation*}
Since $\mathcal{I}^{(t,s,a,s')}$ contains $n^{(t,s,a,s')}$ elements, the submodular polyhedron $\mathcal{P}^{(t,s,a,s')}$ contains $2^{n^{(t,s,a,s')}}$ submodular inequalities. Let $n_y=\sum_{t\in\mathbb{T}}\sum_{s\in\mathbb{S}}
|\mathcal{A}_t(s)|+|\mathbb{S}|$ denote the total number of dual
variables. Let $\mathbf{y}\in\mathbb{R}^{n_y}$ denote the decision
vector of dimension $n_y$ formed by all dual variables $y_t(s,a)$ and
$y_T(s)$. We denote the feasible region of the formulation \eqref{MDP_dual_LP} by $\mathcal{Y}$, where $\mathbf{y}\in\mathcal{Y}$ collects all dual variables. The dual formulation \eqref{MDP_dual_LP} provides another view of our MDP model \eqref{MDP_model}. The local optimal revenues appear only as objective coefficients, and the feasible region $\mathcal{Y}$ does not depend on these revenues. The following lemma provides properties of the set $\mathcal{Y}$ in the formulation \eqref{MDP_dual_LP}. We use this lemma in later proofs.

\begin{lemma}
\label{lem:dual_probability_flow_bounds}
Based on the formulation \eqref{MDP_dual_LP}, for every $\mathbf{y}\in\mathcal{Y}$, we have $\sum_{s\in\mathbb{S}}
\sum_{a\in\mathcal{A}_t(s)}
y_t(s,a)
=
1$ for all $t\in\mathbb{T}$ and $\sum_{s\in\mathbb{S}}
y_T(s)
=
1$.
Thus, we have $0
\leq
y_t(s,a)
\leq
1$ for all $a\in\mathcal{A}_t(s),\
s\in\mathbb{S},\
t\in\mathbb{T}$, and $0
\leq
y_T(s)
\leq
1$ for all $s\in\mathbb{S}$.
\end{lemma}

\begin{proof}
We prove the result by mathematical induction. Based on the formulation \eqref{MDP_dual_LP}, we consider an arbitrary feasible solution
$\bar{\mathbf{y}}\in\mathcal{Y}$. The initial 
constraints give $\sum_{s\in\mathbb{S}}
\sum_{a\in\mathcal{A}_0(s)}
\bar y_0(s,a)
=
1$.
For every $t\in\mathbb{T}\setminus\{0\}$, suppose that
$\sum_{s\in\mathbb{S}}
\sum_{a\in\mathcal{A}_{t-1}(s)}
\bar y_{t-1}(s,a)=1$.
The summation of the constraints at epoch $t$ over
$s\in\mathbb{S}$ gives
\begin{equation*}
\sum_{s\in\mathbb{S}}
\sum_{a\in\mathcal{A}_t(s)}
\bar y_t(s,a)
=
\sum_{\widetilde{s}\in\mathbb{S}}
\sum_{\widetilde{a}\in\mathcal{A}_{t-1}(\widetilde{s})}
\left(
\sum_{s\in\mathbb{S}}
P_{t-1}(s\mid\widetilde{s},\widetilde{a})
\right)
\bar y_{t-1}(\widetilde{s},\widetilde{a})
=
\sum_{\widetilde{s}\in\mathbb{S}}
\sum_{\widetilde{a}\in\mathcal{A}_{t-1}(\widetilde{s})}
\bar y_{t-1}(\widetilde{s},\widetilde{a})
=
1.
\end{equation*}
The mathematical induction over $t$ therefore gives $\sum_{s\in\mathbb{S}}
\sum_{a\in\mathcal{A}_t(s)}
\bar y_t(s,a)
=
1$ for all $t\in\mathbb{T}$. The summation of the terminal 
constraints over $s\in\mathbb{S}$ gives
\begin{equation*}
\sum_{s\in\mathbb{S}}
\bar y_T(s)
=
\sum_{\widetilde{s}\in\mathbb{S}}
\sum_{\widetilde{a}\in\mathcal{A}_{T-1}(\widetilde{s})}
\left(
\sum_{s\in\mathbb{S}}
P_{T-1}(s\mid\widetilde{s},\widetilde{a})
\right)
\bar y_{T-1}(\widetilde{s},\widetilde{a})
=
\sum_{\widetilde{s}\in\mathbb{S}}
\sum_{\widetilde{a}\in\mathcal{A}_{T-1}(\widetilde{s})}
\bar y_{T-1}(\widetilde{s},\widetilde{a})
=
1.
\end{equation*}
The nonnegativity constraints of the formulation
\eqref{MDP_dual_LP} give
$\bar y_t(s,a)\geq0$ for all
$a\in\mathcal{A}_t(s)$, $s\in\mathbb{S}$, and
$t\in\mathbb{T}$. The formulation \eqref{MDP_dual_LP} declares
$\bar y_T(s)\in\mathbb{R}$; however, the terminal 
constraints give $\bar y_T(s)
=
\sum_{\widetilde{s}\in\mathbb{S}}
\sum_{\widetilde{a}\in\mathcal{A}_{T-1}(\widetilde{s})}
P_{T-1}(s\mid\widetilde{s},\widetilde{a})
\bar y_{T-1}(\widetilde{s},\widetilde{a})
\geq
0$ for all $s\in\mathbb{S}$,
where the inequality follows from
$P_{T-1}(s\mid\widetilde{s},\widetilde{a})\geq0$ and
$\bar y_{T-1}(\widetilde{s},\widetilde{a})\geq0$. Since every $\bar y_t(s,a)$ and $\bar y_T(s)$ is nonnegative and
each associated collection sums to one, every variable is bounded
above by one. Given any $t\in\mathbb{T}$, $s\in\mathbb{S}$, and
$a\in\mathcal{A}_t(s)$, we have $0
\leq
\bar y_t(s,a)
\leq
\sum_{\widetilde{s}\in\mathbb{S}}
\sum_{\widetilde{a}\in\mathcal{A}_t(\widetilde{s})}
\bar y_t(\widetilde{s},\widetilde{a})
=
1$. Given any $s\in\mathbb{S}$, we also have $0
\leq
\bar y_T(s)
\leq
\sum_{\widetilde{s}\in\mathbb{S}}
\bar y_T(\widetilde{s})
=
1$. Since
$\bar{\mathbf{y}}\in\mathcal{Y}$ is arbitrary, the above bounds
hold for every $\mathbf{y}\in\mathcal{Y}$. This completes the
proof.
\end{proof}

Next, we integrate the local revenue problems into the dual formulation \eqref{MDP_dual_LP} using the submodular inequalities \eqref{submodular_inequality}. This gives a mixed-integer nonlinear program (MINLP). Let $Z_{\mathrm{D}}$ denote the optimal objective value 
of the formulation \eqref{MDP_dual_LP}. The following proposition shows that the MINLP and the dual LP have the same optimal value $Z_{\mathrm{D}}$.

\begin{proposition}\label{prop:exact_integrated_dual_formulation}
Let $Z_{\mathrm{MINLP}}$
denote the optimal objective value of
\begin{equation}\label{prop_integrated_dual_MINLP}
\begin{aligned}
\max \quad
&\sum_{t\in\mathbb{T}}\sum_{s\in\mathbb{S}}\sum_{a\in\mathcal{A}_t(s)}
\sum_{s'\in\mathbb{S}}
P_t(s'\mid s,a)\theta^{(t,s,a,s')}y_t(s,a)
+\sum_{s\in\mathbb{S}}r_T(s)y_T(s)\\
\text{s.t.}\quad
&\mathbf{y} \in\mathcal{Y},\\
&(\mathbf{x}^{(t,s,a,s')},\theta^{(t,s,a,s')})
\in\mathcal{P}^{(t,s,a,s')}
&&\forall a\in\mathcal{A}_t(s),\ s,s'\in\mathbb{S},\ t\in\mathbb{T},\\
&\mathbf{x}^{(t,s,a,s')}
\in
\mathcal{X}^{(t,s,a,s')}
\cap
\mathbb{B}^{n^{(t,s,a,s')}}
&&\forall a\in\mathcal{A}_t(s),\ s,s'\in\mathbb{S},\ t\in\mathbb{T}.
\end{aligned}
\end{equation}
Then, the formulation \eqref{prop_integrated_dual_MINLP} is a mixed-integer bilinear formulation with $Z_{\mathrm{MINLP}}=Z_{\mathrm{D}}=V_0^*(s_0)$.
\end{proposition}

\begin{proof}
The only nonlinear terms of the formulation
\eqref{prop_integrated_dual_MINLP} are the products
$\theta^{(t,s,a,s')}y_t(s,a)$, where both variables are continuous.
Since each local vector $\mathbf{x}^{(t,s,a,s')}$ is binary, the
formulation \eqref{prop_integrated_dual_MINLP} has a mixed-integer
bilinear structure.

We first show that $Z_{\mathrm{MINLP}}\leq Z_{\mathrm{D}}$. We consider
 arbitrary feasible solutions $\bar{\mathbf{y}}$,
$\bar{\mathbf{x}}^{(t,s,a,s')}$, and
$\bar{\theta}^{(t,s,a,s')}$ to the formulation
\eqref{prop_integrated_dual_MINLP}. The feasibility of the selected
solution gives $\big(
\bar{\mathbf{x}}^{(t,s,a,s')},
\bar{\theta}^{(t,s,a,s')}
\big)
\in
\mathcal{P}^{(t,s,a,s')}$ and  $\bar{\mathbf{x}}^{(t,s,a,s')}
\in
\mathcal{X}^{(t,s,a,s')}
\cap
\mathbb{B}^{n^{(t,s,a,s')}}$.
The submodular inequality \eqref{submodular_inequality} indexed by
$\overline{X}^{(t,s,a,s')}$ provides
\begin{equation*}
\begin{aligned}
\bar{\theta}^{(t,s,a,s')}
\leq
f^{(t,s,a,s')}(\overline{X}^{(t,s,a,s')})
&-
\underbrace{
\sum_{i\in\overline{X}^{(t,s,a,s')}}
\rho_i^{(t,s,a,s')}
(\mathcal{I}^{(t,s,a,s')}\setminus\{i\})
(1-\bar{x}_i^{(t,s,a,s')})
}_{=\,0}\\
&+
\underbrace{
\sum_{i\in\mathcal{I}^{(t,s,a,s')}\setminus
\overline{X}^{(t,s,a,s')}}
\rho_i^{(t,s,a,s')}(\overline{X}^{(t,s,a,s')})
\bar{x}_i^{(t,s,a,s')}
}_{=\,0}\\
=\;&
f^{(t,s,a,s')}(\bar{\mathbf{x}}^{(t,s,a,s')})
\leq
r_t(s,a,s').
\end{aligned}
\end{equation*}
Therefore, the local variables of the selected feasible solution
satisfy
\begin{equation}\label{prop_revenue_upper_bound}
\bar{\theta}^{(t,s,a,s')}
\leq
f^{(t,s,a,s')}
(\bar{\mathbf{x}}^{(t,s,a,s')})
\leq
r_t(s,a,s')
\quad
\forall a\in\mathcal{A}_t(s),\
s,s'\in\mathbb{S},\
t\in\mathbb{T},
\end{equation}
where the second inequality follows from the feasibility of
$\bar{\mathbf{x}}^{(t,s,a,s')}$ for the local problem
\eqref{MDP_transition_revenue}. From Lemma \ref{lem:dual_probability_flow_bounds}, since
$P_t(s'\mid s,a)\bar{y}_t(s,a)\geq 0$, we multiply the relation
\eqref{prop_revenue_upper_bound} by
$P_t(s'\mid s,a)\bar{y}_t(s,a)$ and sum over all transition indices.
Adding the terminal revenue term to both sides gives
\begin{equation}\label{prop_dual_objective_bound}
\begin{aligned}
&
\sum_{t\in\mathbb{T}}
\sum_{s\in\mathbb{S}}
\sum_{a\in\mathcal{A}_t(s)}
\sum_{s'\in\mathbb{S}}
P_t(s'\mid s,a)
\bar{\theta}^{(t,s,a,s')}
\bar{y}_t(s,a)
+
\sum_{s\in\mathbb{S}}
r_T(s)\bar{y}_T(s)
\\
\leq\;&
\sum_{t\in\mathbb{T}}
\sum_{s\in\mathbb{S}}
\sum_{a\in\mathcal{A}_t(s)}
\sum_{s'\in\mathbb{S}}
P_t(s'\mid s,a)
r_t(s,a,s')
\bar{y}_t(s,a)
+
\sum_{s\in\mathbb{S}}
r_T(s)\bar{y}_T(s)
\leq
Z_{\mathrm{D}}.
\end{aligned}
\end{equation}
The first inequality follows from the relation
\eqref{prop_revenue_upper_bound}. The second inequality follows from
$\bar{\mathbf{y}}\in\mathcal{Y}$ and the definition of
$Z_{\mathrm{D}}$. Since the selected feasible solutions to the formulation 
\eqref{prop_integrated_dual_MINLP} are arbitrary, the relation
\eqref{prop_dual_objective_bound} holds for all feasible solutions.
Therefore,
$Z_{\mathrm{MINLP}}\leq Z_{\mathrm{D}}$.

We next show that the upper bound $Z_{\mathrm{D}}$ is attained by a
feasible solution to the formulation
\eqref{prop_integrated_dual_MINLP}.  Since $\mathbb{S}$ is finite, the formulation
\eqref{MDP_dual_LP} is an LP with finite variables and constraints. Its feasible region
is nonempty, and its optimal objective value $Z_{\mathrm{D}}$ is
finite. Therefore, the formulation \eqref{MDP_dual_LP} admits an
optimal solution. Let $\mathbf{y}^*\in\mathcal{Y}$ denote an
optimal solution to the
formulation \eqref{MDP_dual_LP}. We consider an
optimal solution $\mathbf{x}^{*(t,s,a,s')}$ to each local submodular maximization problem
\eqref{MDP_transition_revenue}. We define
\begin{equation}\label{prop_optimal_revenue_selection}
\theta^{*(t,s,a,s')}
=
f^{(t,s,a,s')}
(\mathbf{x}^{*(t,s,a,s')})
=
r_t(s,a,s')
\quad
\forall a\in\mathcal{A}_t(s),\
s,s'\in\mathbb{S},\
t\in\mathbb{T}.
\end{equation}
Based on the relation
\eqref{prop_optimal_revenue_selection}, the objective value of the
formulation \eqref{prop_integrated_dual_MINLP} is
\begin{equation*}\label{prop_MINLP_attains_dual_value}
\begin{aligned}
&
\sum_{t\in\mathbb{T}}
\sum_{s\in\mathbb{S}}
\sum_{a\in\mathcal{A}_t(s)}
\sum_{s'\in\mathbb{S}}
P_t(s'\mid s,a)
\theta^{*(t,s,a,s')}y_t^*(s,a)
+
\sum_{s\in\mathbb{S}}
r_T(s)y_T^*(s)
\\
=\;&
\sum_{t\in\mathbb{T}}
\sum_{s\in\mathbb{S}}
\sum_{a\in\mathcal{A}_t(s)}
\sum_{s'\in\mathbb{S}}
P_t(s'\mid s,a)
r_t(s,a,s')y_t^*(s,a)
+
\sum_{s\in\mathbb{S}}
r_T(s)y_T^*(s)
=
Z_{\mathrm{D}}.
\end{aligned}
\end{equation*}
Thus, the formulation
\eqref{prop_integrated_dual_MINLP} has a feasible solution with
objective value $Z_{\mathrm{D}}$. Together with the previous
 upper bound $Z_{\mathrm{MINLP}}
\le
Z_{\mathrm{D}}$, this gives
$Z_{\mathrm{MINLP}}
=
Z_{\mathrm{D}}.$
Also, the strong LP duality theorem gives
$Z_{\mathrm{D}}
=
V_0^*(s_0).$
Therefore, the final relation
$Z_{\mathrm{MINLP}}
=
Z_{\mathrm{D}}
=
V_0^*(s_0)$ holds.
This completes the proof.
\end{proof}

Next, we show that the bilinear formulation \eqref{prop_integrated_dual_MINLP} can be replaced by a mixed-integer linear program (MILP) with the same optimal value. This reformulation uses the properties of $\mathcal Y$, the concepts of  McCormick linearization \cite{Adams1993,LuedtkeMP2012,McC76}, and the submodular inequalities \eqref{submodular_inequality}. Here, we also refer the reader to \cite{MINLP2013,KleinertEJCO2021} for standard literature on techniques for MINLP and bilinear optimization. For each $t\in\mathbb{T}$, $s,s'\in\mathbb{S}$, and
$a\in\mathcal{A}_t(s)$, we introduce a continuous variable
$d^{(t,s,a,s')}\in\mathbb{R}$ and a continuous vector
$\mathbf{w}^{(t,s,a,s')}
\in\mathbb{R}^{n^{(t,s,a,s')}}$. The variables
$d^{(t,s,a,s')}$ and $w_i^{(t,s,a,s')}$ are used in the linear
reformulation of the products
$\theta^{(t,s,a,s')}y_t(s,a)$ and
$x_i^{(t,s,a,s')}y_t(s,a)$, respectively. For these indices, we have
the perspective submodular polyhedron
\begin{equation*}
\mathcal{D}^{(t,s,a,s')}
=
\left\{
\begin{aligned}
&\big(
d^{(t,s,a,s')},
y_t(s,a),
\mathbf{w}^{(t,s,a,s')}
\big)
\in
\mathbb{R}\times[0,1]\times
\mathbb{R}^{n^{(t,s,a,s')}}:\\
&d^{(t,s,a,s')}
\leq
f^{(t,s,a,s')}(\overline{X}^{(t,s,a,s')})
y_t(s,a)\\
&\quad-
\sum_{i\in\overline{X}^{(t,s,a,s')}}
\rho_i^{(t,s,a,s')}
(\mathcal{I}^{(t,s,a,s')}\setminus\{i\})
\big(
y_t(s,a)-w_i^{(t,s,a,s')}
\big)\\
&\quad+
\sum_{i\in\mathcal{I}^{(t,s,a,s')}
\setminus\overline{X}^{(t,s,a,s')}}
\rho_i^{(t,s,a,s')}
(\overline{X}^{(t,s,a,s')})
w_i^{(t,s,a,s')}\\
&\hspace{25mm}
\forall\,
\overline{X}^{(t,s,a,s')}
\subseteq
\mathcal{I}^{(t,s,a,s')}
\end{aligned}
\right\}.
\end{equation*}
For the same indices, we have the McCormick polyhedron
\begin{equation*}
\mathcal{W}^{(t,s,a,s')}
=
\left\{
\begin{aligned}
&\big(
\mathbf{x}^{(t,s,a,s')},
y_t(s,a),
\mathbf{w}^{(t,s,a,s')}
\big)
\in
\mathbb{B}^{n^{(t,s,a,s')}}\times[0,1]\times
\mathbb{R}^{n^{(t,s,a,s')}}:\\
&0\leq w_i^{(t,s,a,s')}\leq y_t(s,a),\\
&w_i^{(t,s,a,s')}
\leq
x_i^{(t,s,a,s')},\\
&w_i^{(t,s,a,s')}
\geq
y_t(s,a)+x_i^{(t,s,a,s')}-1
\quad
\forall i\in\mathcal{I}^{(t,s,a,s')}
\end{aligned}
\right\}.
\end{equation*}
In the following proposition, each local constraint is applied for all $t\in\mathbb{T}$, $s,s'\in\mathbb{S}$, and $a\in\mathcal{A}_t(s)$.

\begin{proposition}
\label{prop:equivalent_MILP_reformulation}
Let $Z_{\mathrm{MILP}}$ denote the optimal
objective value of
\begin{equation}\label{MDP_submodular_dual_MILP}
\begin{aligned}
\max\quad
&
\sum_{t\in\mathbb{T}}
\sum_{s\in\mathbb{S}}
\sum_{a\in\mathcal{A}_t(s)}
\sum_{s'\in\mathbb{S}}
P_t(s'\mid s,a)d^{(t,s,a,s')}
+
\sum_{s\in\mathbb{S}}r_T(s)y_T(s)
\\
\text{s.t.}\quad
&\mathbf{y}\in\mathcal{Y},\\
&\mathbf{x}^{(t,s,a,s')}
\in
\mathcal{X}^{(t,s,a,s')}
\cap
\mathbb{B}^{n^{(t,s,a,s')}},\\
&\big(
d^{(t,s,a,s')},
y_t(s,a),
\mathbf{w}^{(t,s,a,s')}
\big)
\in
\mathcal{D}^{(t,s,a,s')},\\
&\big(
\mathbf{x}^{(t,s,a,s')},
y_t(s,a),
\mathbf{w}^{(t,s,a,s')}
\big)
\in
\mathcal{W}^{(t,s,a,s')}.
\end{aligned}
\end{equation}
Then, the formulation \eqref{MDP_submodular_dual_MILP} is a
mixed-integer linear reformulation of the formulation
\eqref{prop_integrated_dual_MINLP} with $Z_{\mathrm{MILP}}
=
Z_{\mathrm{MINLP}}
=
Z_{\mathrm{D}}
=
V_0^*(s_0)$.
\end{proposition}

\begin{proof}
The sets $\mathcal{Y}$, $\mathcal{X}^{(t,s,a,s')}$,
$\mathcal{D}^{(t,s,a,s')}$, and
$\mathcal{W}^{(t,s,a,s')}$ are described by linear constraints in the formulation
\eqref{MDP_submodular_dual_MILP}. In addition, since every local vector
$\mathbf{x}^{(t,s,a,s')}$ is binary, the formulation
\eqref{MDP_submodular_dual_MILP} is a mixed-integer linear program.

We first verify the
relation represented by the polyhedron
$\mathcal{W}^{(t,s,a,s')}$ for binary
$\mathbf{x}^{(t,s,a,s')}$. Based on Lemma \ref{lem:dual_probability_flow_bounds}, for every $\mathbf{y}\in\mathcal{Y}$, the constraints defining
$\mathcal{Y}$ give
$0\leq y_t(s,a)\leq 1$
for all $a\in\mathcal{A}_t(s)$,
$s\in\mathbb{S}$, and $t\in\mathbb{T}$. Given  $t\in\mathbb{T}$,
$s,s'\in\mathbb{S}$, $a\in\mathcal{A}_t(s)$, and
$i\in\mathcal{I}^{(t,s,a,s')}$, if
$x_i^{(t,s,a,s')}=0$, then
$0
\leq
w_i^{(t,s,a,s')}
\leq
x_i^{(t,s,a,s')}
=
0$.
If $x_i^{(t,s,a,s')}=1$, then
$y_t(s,a)
=
y_t(s,a)+x_i^{(t,s,a,s')}-1
\leq
w_i^{(t,s,a,s')}
\leq
y_t(s,a)$.
Therefore, every feasible solution in
$\mathcal{W}^{(t,s,a,s')}$ with binary
$\mathbf{x}^{(t,s,a,s')}$ satisfies
$w_i^{(t,s,a,s')}
=
x_i^{(t,s,a,s')}y_t(s,a)$
for all $i\in\mathcal{I}^{(t,s,a,s')}$. On the other hand, suppose that
$w_i^{(t,s,a,s')}=
x_i^{(t,s,a,s')}y_t(s,a)$. Based on
$x_i^{(t,s,a,s')}\in\mathbb{B}$ and
$y_t(s,a)\in[0,1]$, we have
\begin{equation*}
\begin{aligned}
&0
\leq
w_i^{(t,s,a,s')}
=
x_i^{(t,s,a,s')}y_t(s,a)
\leq
y_t(s,a),\\
&w_i^{(t,s,a,s')}
=
x_i^{(t,s,a,s')}y_t(s,a)
\leq
x_i^{(t,s,a,s')},\\
&w_i^{(t,s,a,s')}
-
y_t(s,a)
-
x_i^{(t,s,a,s')}
+
1
=
\big(1-x_i^{(t,s,a,s')}\big)
\big(1-y_t(s,a)\big)
\geq
0,
\end{aligned}
\end{equation*}
where the last relation provides the validity of the constraint
$w_i^{(t,s,a,s')}
\geq
y_t(s,a)+x_i^{(t,s,a,s')}-1$
of $\mathcal{W}^{(t,s,a,s')}$. Therefore, the polyhedron
$\mathcal{W}^{(t,s,a,s')}$ exactly represents the products
$x_i^{(t,s,a,s')}y_t(s,a)$ for binary
$\mathbf{x}^{(t,s,a,s')}$ and $y_t(s,a)\in[0,1]$. In other
words, for
$\mathbf{x}^{(t,s,a,s')}
\in\mathbb{B}^{n^{(t,s,a,s')}}$
and $y_t(s,a)\in[0,1]$, we have the relation
\begin{equation}\label{prop_McCormick_exactness}
\big(
\mathbf{x}^{(t,s,a,s')},
y_t(s,a),
\mathbf{w}^{(t,s,a,s')}
\big)
\in
\mathcal{W}^{(t,s,a,s')}
 \text{ if and only if } 
w_i^{(t,s,a,s')}
=
x_i^{(t,s,a,s')}y_t(s,a)
\quad
\forall i\in\mathcal{I}^{(t,s,a,s')}.
\end{equation}

We next verify that the product substitutions define a feasible point
in $\mathcal{D}^{(t,s,a,s')}$. Suppose that
$\big(
\mathbf{x}^{(t,s,a,s')},
\theta^{(t,s,a,s')}
\big)
\in
\mathcal{P}^{(t,s,a,s')}$
and $y_t(s,a)\in[0,1]$. Define
$d^{(t,s,a,s')}=
\theta^{(t,s,a,s')}y_t(s,a)$ and
$w_i^{(t,s,a,s')}=
x_i^{(t,s,a,s')}y_t(s,a)$ for every
$i\in\mathcal{I}^{(t,s,a,s')}$. These definitions give
\begin{equation*}
\begin{aligned}
\theta^{(t,s,a,s')}y_t(s,a)
&=
d^{(t,s,a,s')},\\
x_i^{(t,s,a,s')}y_t(s,a)
&=
w_i^{(t,s,a,s')},\\
\big(1-x_i^{(t,s,a,s')}\big)y_t(s,a)
&=
y_t(s,a)-w_i^{(t,s,a,s')}.
\end{aligned}
\end{equation*}
We multiply every submodular inequality in \eqref{submodular_inequality} by $y_t(s,a)\geq 0$. We have a relation that if  
\begin{equation}\label{prop_perspective_lift}
\left.
\begin{aligned}
&
\big(
\mathbf{x}^{(t,s,a,s')},
\theta^{(t,s,a,s')}
\big)
\in
\mathcal{P}^{(t,s,a,s')},\\
&y_t(s,a)\in[0,1],\\
&d^{(t,s,a,s')}
=
\theta^{(t,s,a,s')}y_t(s,a),\\
&w_i^{(t,s,a,s')}
=
x_i^{(t,s,a,s')}y_t(s,a)
\quad
\forall i\in\mathcal{I}^{(t,s,a,s')}
\end{aligned}
\right\}
\text{hold, then }
\big(
d^{(t,s,a,s')},
y_t(s,a),
\mathbf{w}^{(t,s,a,s')}
\big)
\in
\mathcal{D}^{(t,s,a,s')}.
\end{equation}
The formulation \eqref{MDP_submodular_dual_MILP} contains the
membership constraint in
$\mathcal{D}^{(t,s,a,s')}$, but it contains neither
$\theta^{(t,s,a,s')}$ nor the equality
$d^{(t,s,a,s')}=
\theta^{(t,s,a,s')}y_t(s,a)$.

The above analysis verifies that the product substitutions define
feasible points in $\mathcal{W}^{(t,s,a,s')}$ and
$\mathcal{D}^{(t,s,a,s')}$. We next provide an upper bound on the
objective value of every feasible solution to the formulation
\eqref{MDP_submodular_dual_MILP}. Specifically, we show that
$Z_{\mathrm{MILP}}\leq Z_{\mathrm{MINLP}}$. We consider any feasible solutions $\bar{\mathbf{y}}$,
$\bar{\mathbf{x}}^{(t,s,a,s')}$,
$\bar{d}^{(t,s,a,s')}$, and
$\bar{\mathbf{w}}^{(t,s,a,s')}$ to the formulation
\eqref{MDP_submodular_dual_MILP}. Based on the relation
\eqref{prop_McCormick_exactness}, we have
$\bar{w}_i^{(t,s,a,s')}
=
\bar{y}_t(s,a)\bar{x}_i^{(t,s,a,s')}$
for all
$i\in\mathcal{I}^{(t,s,a,s')}$. We select the inequality in
$\mathcal{D}^{(t,s,a,s')}$ indexed by
$\overline{X}^{(t,s,a,s')}$. The corresponding substitution gives
\begin{equation}\label{prop_MILP_local_revenue_bound}
\begin{aligned}
\bar{d}^{(t,s,a,s')}
\leq\;&
f^{(t,s,a,s')}(\overline{X}^{(t,s,a,s')})\bar{y}_t(s,a)\\
&-
\sum_{i\in\overline{X}^{(t,s,a,s')}}
\rho_i^{(t,s,a,s')}
\big(\mathcal{I}^{(t,s,a,s')}\setminus\{i\}\big)
\left(
\bar{y}_t(s,a)
-
\underbrace{
\bar{y}_t(s,a)\bar{x}_i^{(t,s,a,s')}
}_{=\,\bar{w}_i^{(t,s,a,s')}}
\right)\\
&+
\sum_{i\in\mathcal{I}^{(t,s,a,s')}
\setminus\overline{X}^{(t,s,a,s')}}
\rho_i^{(t,s,a,s')}
(\overline{X}^{(t,s,a,s')})
\underbrace{
\bar{y}_t(s,a)\bar{x}_i^{(t,s,a,s')}
}_{=\,\bar{w}_i^{(t,s,a,s')}}\\
=\;&
f^{(t,s,a,s')}(\overline{X}^{(t,s,a,s')})\bar{y}_t(s,a)\\
&-
\underbrace{
\sum_{i\in\overline{X}^{(t,s,a,s')}}
\rho_i^{(t,s,a,s')}
\big(\mathcal{I}^{(t,s,a,s')}\setminus\{i\}\big)
\big(\bar{y}_t(s,a)-\bar{y}_t(s,a)\big)
}_{=\,0}\\
&+
\underbrace{
\sum_{i\in\mathcal{I}^{(t,s,a,s')}
\setminus\overline{X}^{(t,s,a,s')}}
\rho_i^{(t,s,a,s')}
(\overline{X}^{(t,s,a,s')})\,0
}_{=\,0}\\
=\;&
f^{(t,s,a,s')}(\overline{X}^{(t,s,a,s')})\bar{y}_t(s,a)
=
f^{(t,s,a,s')}
(\bar{\mathbf{x}}^{(t,s,a,s')})\bar{y}_t(s,a).
\end{aligned}
\end{equation}
We define $\bar{\theta}^{(t,s,a,s')}$ as the largest feasible value
\begin{equation}\label{prop_MILP_to_MINLP_mapping}
\begin{aligned}
\bar{\theta}^{(t,s,a,s')}
=
\max
\left\{
\theta^{(t,s,a,s')}\in\mathbb{R}:
\big(
\bar{\mathbf{x}}^{(t,s,a,s')},
\theta^{(t,s,a,s')}
\big)
\in
\mathcal{P}^{(t,s,a,s')}
\right\}
=
f^{(t,s,a,s')}
(\bar{\mathbf{x}}^{(t,s,a,s')}),
\end{aligned}
\end{equation}
where
$\big(
\bar{\mathbf{x}}^{(t,s,a,s')},
\bar{\theta}^{(t,s,a,s')}
\big)
\in
\mathcal{P}^{(t,s,a,s')}$. Based on the relations
\eqref{prop_MILP_local_revenue_bound} and
\eqref{prop_MILP_to_MINLP_mapping}, we have
\begin{equation}\label{prop_MILP_product_bound}
\bar{d}^{(t,s,a,s')}
\leq
f^{(t,s,a,s')}
(\bar{\mathbf{x}}^{(t,s,a,s')})\bar{y}_t(s,a)
=
\bar{\theta}^{(t,s,a,s')}\bar{y}_t(s,a).
\end{equation}
Since
$P_t(s'\mid s,a)\geq 0$, multiplying the relation 
\eqref{prop_MILP_product_bound} by $P_t(s'\mid s,a)$ gives
$P_t(s'\mid s,a)\bar{d}^{(t,s,a,s')}
\leq
P_t(s'\mid s,a)
\bar{\theta}^{(t,s,a,s')}\bar{y}_t(s,a)$.
Summing these inequalities over all transitions and adding the
terminal revenue term gives
\begin{equation}\label{prop_MILP_objective_bound}
\begin{aligned}
&
\sum_{t\in\mathbb{T}}
\sum_{s\in\mathbb{S}}
\sum_{a\in\mathcal{A}_t(s)}
\sum_{s'\in\mathbb{S}}
P_t(s'\mid s,a)\bar{d}^{(t,s,a,s')}
+
\sum_{s\in\mathbb{S}}r_T(s)\bar{y}_T(s)
\\
\leq\;&
\sum_{t\in\mathbb{T}}
\sum_{s\in\mathbb{S}}
\sum_{a\in\mathcal{A}_t(s)}
\sum_{s'\in\mathbb{S}}
P_t(s'\mid s,a)
\bar{\theta}^{(t,s,a,s')}\bar{y}_t(s,a)
+
\sum_{s\in\mathbb{S}}r_T(s)\bar{y}_T(s)
\leq
Z_{\mathrm{MINLP}},
\end{aligned}
\end{equation}
where the last inequality follows from the feasibility of the
constructed solution to the formulation
\eqref{prop_integrated_dual_MINLP}. Since the relation
\eqref{prop_MILP_objective_bound} holds for every feasible solution
to the formulation \eqref{MDP_submodular_dual_MILP}, we obtain
$Z_{\mathrm{MILP}}
\leq
Z_{\mathrm{MINLP}}$.

We next show that the value $Z_{\mathrm{MINLP}}$ is attained by a
feasible solution to the formulation
\eqref{MDP_submodular_dual_MILP}. Proposition
\ref{prop:exact_integrated_dual_formulation} gives the existence of
an optimal solution to the formulation
\eqref{prop_integrated_dual_MINLP}. Let
$\mathbf{y}^*$,
$\mathbf{x}^{*(t,s,a,s')}$, and
$\theta^{*(t,s,a,s')}$ denote such an optimal solution. We define
$d^{*(t,s,a,s')}
=
\theta^{*(t,s,a,s')}y_t^*(s,a)$ and
$w_i^{*(t,s,a,s')}
=
x_i^{*(t,s,a,s')}y_t^*(s,a)$.
Based on the
definition
$d^{*(t,s,a,s')}
=
\theta^{*(t,s,a,s')}y_t^*(s,a)$, its objective value satisfies
\begin{equation*}\label{prop_MILP_attains_MINLP_value}
\begin{aligned}
&
\sum_{t\in\mathbb{T}}
\sum_{s\in\mathbb{S}}
\sum_{a\in\mathcal{A}_t(s)}
\sum_{s'\in\mathbb{S}}
P_t(s'\mid s,a)d^{*(t,s,a,s')}
+
\sum_{s\in\mathbb{S}}r_T(s)y_T^*(s)
\\
=\;&
\sum_{t\in\mathbb{T}}
\sum_{s\in\mathbb{S}}
\sum_{a\in\mathcal{A}_t(s)}
\sum_{s'\in\mathbb{S}}
P_t(s'\mid s,a)
\theta^{*(t,s,a,s')}y_t^*(s,a)
+
\sum_{s\in\mathbb{S}}r_T(s)y_T^*(s)
=
Z_{\mathrm{MINLP}}.
\end{aligned}
\end{equation*}
Thus, the formulation \eqref{MDP_submodular_dual_MILP} has a feasible solution with
objective value $Z_{\mathrm{MINLP}}$. Together with
$Z_{\mathrm{MILP}}
\leq
Z_{\mathrm{MINLP}}$, this gives
$Z_{\mathrm{MILP}}
=
Z_{\mathrm{MINLP}}$.
Proposition \ref{prop:exact_integrated_dual_formulation} gives
$Z_{\mathrm{MINLP}}
=
Z_{\mathrm{D}}
=
V_0^*(s_0)$.
This completes the proof.
\end{proof}
Based on the dual formulation \eqref{MDP_dual_LP} and  Proposition~\ref{prop:equivalent_MILP_reformulation}, we obtain an exact MILP formulation \eqref{MDP_submodular_dual_MILP} of our MDP model \eqref{MDP_model}, providing a basis for further study using mixed-integer programming techniques. That is, for each local revenue problem \eqref{MDP_transition_revenue}, the formulation introduces auxiliary variables and submodular inequalities for all subsets of its ground set. The formulation makes the local constraints explicit and links the local problems \eqref{MDP_transition_revenue} through the MDP dual variables. This structure may suggest further research topics for tighter polyhedral formulations and decomposition methods. Here, we refer to \cite{SMIP2017,NW88,SenHigle1999} for background on these topics.  We note that the local optimal revenues enter the dual objective as coefficients. In this study, we also use the dual formulation \eqref{MDP_dual_LP} to analyze the sensitivity of the MDP optimal value to changes in the local revenues. The next section uses incumbent revenues and valid upper bounds for the local problems to evaluate the effect of local optimality gaps on the MDP optimal value without requiring all local problems to be solved to optimality.

\section{A Dual Sensitivity Analysis of the MDP Optimality Gap}
The sensitivity analysis of this section is motivated by classical discussions in optimization textbooks. For example, Chapter II.3 of \cite{NW88} discusses duality and the value function, Chapter 5 of \cite{BT1997} discusses sensitivity analysis in linear programming, and Section 10 of \cite{Schrijver98} uses dual solutions to
analyze changes in the LP optimal value. These discussions motivate our analysis of
changes in the MDP optimal value. We use local incumbent revenues and valid upper bounds in this analysis. In our setting, the local revenues appear only in the objective function and do not affect the feasible region. This dual structure provides a natural starting point for the sensitivity analysis in this study. Recall that each local revenue problem \eqref{MDP_transition_revenue} is NP-hard.
If a method provides an incumbent revenue and a valid optimality gap for each local revenue problem \eqref{MDP_transition_revenue}, we can use this information and the dual formulation \eqref{MDP_dual_LP} to provide a dual sensitivity analysis of the MDP optimality gap of the model \eqref{MDP_model}. 

Formally, based on the dual formulation \eqref{MDP_dual_LP}, the local
transition revenues appear only as objective coefficients, while 
the set of constraints $\mathcal{Y}$ remains unchanged.
For each local problem \eqref{MDP_transition_revenue}, we assume that a method provides an incumbent local solution 
$\bar{\mathbf{x}}^{(t,s,a,s')}$ and an incumbent local revenue as $\bar{r}_t(s,a,s')
=
f^{(t,s,a,s')}
\big(
\bar{\mathbf{x}}^{(t,s,a,s')}
\big)$. Since $\bar{\mathbf{x}}^{(t,s,a,s')}$ is feasible for the local problem
\eqref{MDP_transition_revenue}, we have
$\bar{r}_t(s,a,s')\leq r_t(s,a,s')$ for the associated transition revenue. We replace the optimal transition revenues $r_t(s,a,s')$ in the dual formulation
\eqref{MDP_dual_LP} with the incumbent local revenues $\bar{r}_t(s,a,s')$. Then, we define $\bar{Z}_{\mathrm{D}}$ as the optimal objective value of the
resulting dual formulation
\begin{equation*}
\bar{Z}_{\mathrm{D}}
=
\max_{\mathbf{y} \in \mathcal{Y}}\quad
\sum_{t\in\mathbb{T}}
\sum_{s\in\mathbb{S}}
\sum_{a\in\mathcal{A}_t(s)}
\sum_{s'\in\mathbb{S}}
P_t(s'\mid s,a)
\bar{r}_t(s,a,s')
y_t(s,a)
+
\sum_{s\in\mathbb{S}}
r_T(s)y_T(s).
\end{equation*}
The dual formulations \eqref{MDP_dual_LP} and
$\bar{Z}_{\mathrm{D}}$ share the feasible
region $\mathcal{Y}$ and the same terminal revenue coefficients.
Since $y_t(s,a)\geq 0$ and
$\bar{r}_t(s,a,s')\leq r_t(s,a,s')$, the comparison of the coefficients provides the baseline optimality bound
$\bar{Z}_{\mathrm{D}}\leq Z_{\mathrm{D}}$. For each $t\in\mathbb{T}$, $s,s'\in\mathbb{S}$, and
$a\in\mathcal{A}_t(s)$, we define $U_t(s,a,s')$ as a valid upper bound
on the optimal objective value of the corresponding local  problem \eqref{MDP_transition_revenue}, where $\bar{r}_t(s,a,s')
\leq
r_t(s,a,s')
\leq
U_t(s,a,s')$.
We define the certified relative optimality gap of each local problem \eqref{MDP_transition_revenue}
based on its incumbent revenue and upper bound as
\begin{equation*}
\delta_t(s,a,s')
=
\begin{cases}
\dfrac{
U_t(s,a,s')
-
\bar r_t(s,a,s')
}{
U_t(s,a,s')
},
&
U_t(s,a,s')>0,\
\bar r_t(s,a,s')\geq0,\\[4mm]
\dfrac{
U_t(s,a,s')
-
\bar r_t(s,a,s')
}{
\left|
U_t(s,a,s')
\right|
},
&
U_t(s,a,s')<0,\\[4mm]
0,
&
U_t(s,a,s')
=
\bar r_t(s,a,s')
=
0,\\[2mm]
+\infty,
&
\bar r_t(s,a,s')<0
\leq
U_t(s,a,s').
\end{cases}
\end{equation*}
If the optimality gap $\delta_t(s,a,s')$ is bounded, the definition above gives $0
\leq
r_t(s,a,s')
-
\bar r_t(s,a,s')
\leq
\delta_t(s,a,s')
\left|
r_t(s,a,s')
\right|.$
Starting with the baseline bound
$\bar Z_{\mathrm{D}}\leq Z_{\mathrm{D}}$, we use the dual
feasible region $\mathcal Y$ and the local optimality gap bounds above to obtain
the following dual optimality gap bounds. The following proposition analyzes the sensitivity of the dual optimal
value to changes in the local revenue coefficients. If each local revenue coefficient is changed to an incumbent revenue with an optimality gap at most $\delta$, then the change in the dual optimal value satisfies
$0\leq Z_{\mathrm{D}}-\bar Z_{\mathrm{D}}
\leq\delta Z_{\mathrm{D}}$.

\begin{proposition}
\label{prop:dual_uniform_local_optimality_gap}
Suppose that there exists $\delta\in(0,1)$ such that $0
\leq
\delta_t(s,a,s')
\leq
\delta$ for all $a\in\mathcal{A}_t(s),\
s,s'\in\mathbb{S},\
t\in\mathbb{T}.$
If the transition and terminal revenues are nonnegative with $r_t(s,a,s')
\geq
0$ 
for all
 $a\in\mathcal{A}_t(s),\
s,s'\in\mathbb{S},\
t\in\mathbb{T}$ and $r_T(s)
\geq
0$
for all
$s\in\mathbb{S}$,
then we have $0
\leq
Z_{\mathrm{D}}
-
\bar Z_{\mathrm{D}}
\leq
\delta Z_{\mathrm{D}}$
\end{proposition}

\begin{proof}
We first establish a baseline relation used throughout
the proof. Consider an arbitrary feasible solution
$\bar{\mathbf{y}}\in\mathcal{Y}$. Based on Lemma
\ref{lem:dual_probability_flow_bounds}, all dual variables $\mathbf{y}\in\mathcal{Y}$ are
nonnegative. Based on the assumption of the local optimality gap condition, we
also have
$r_t(s,a,s')-\bar r_t(s,a,s')\geq0$ for every transition. Let $F(\mathbf{y})$ and $\bar F(\mathbf{y})$ denote the
objective functions of the exact dual formulation and the dual
formulation based on the local incumbent revenues, respectively. That is, since
$P_t(s'\mid s,a)\bar y_t(s,a)\geq0$ from Lemma \ref{lem:dual_probability_flow_bounds}, the difference between the two
dual objective values at $\bar{\mathbf{y}}$ satisfies
\begin{equation*}
F(\bar{\mathbf{y}})
-
\bar F(\bar{\mathbf{y}})
=
\sum_{t\in\mathbb{T}}
\sum_{s\in\mathbb{S}}
\sum_{a\in\mathcal{A}_t(s)}
\sum_{s'\in\mathbb{S}}
P_t(s'\mid s,a)
\left(
r_t(s,a,s')
-
\bar r_t(s,a,s')
\right)
\bar y_t(s,a)
\geq
0.
\end{equation*}
The terminal revenue terms are eliminated because both dual formulations use
the same terminal revenues. Since
$\bar{\mathbf{y}}\in\mathcal{Y}$ is arbitrary, we have
$\bar F(\mathbf{y})\leq F(\mathbf{y})$ for every
$\mathbf{y}\in\mathcal{Y}$. Maximizing both objective functions
over the common feasible region gives the baseline relation $\bar Z_{\mathrm{D}}
=
\max_{\mathbf{y}\in\mathcal{Y}}
\bar F(\mathbf{y})
\leq
\max_{\mathbf{y}\in\mathcal{Y}}
F(\mathbf{y})
=
Z_{\mathrm{D}}$.

Since every transition revenue is
nonnegative, the local optimality gap condition gives $\bar r_t(s,a,s')
\geq
r_t(s,a,s')
-
\delta
\left|
r_t(s,a,s')
\right|
=
\big(
1-\delta
\big)
r_t(s,a,s')$ for all $a\in\mathcal{A}_t(s),\
s,s'\in\mathbb{S},\
t\in\mathbb{T}$. For every $\bar{\mathbf{y}}\in\mathcal{Y}$, the nonnegativity of
the transition probabilities and Lemma
\ref{lem:dual_probability_flow_bounds} give
\begin{equation*}
\begin{aligned}
\bar F(\bar{\mathbf{y}})
\geq\;&
\big(
1-\delta
\big)
\sum_{t\in\mathbb{T}}
\sum_{s\in\mathbb{S}}
\sum_{a\in\mathcal{A}_t(s)}
\sum_{s'\in\mathbb{S}}
P_t(s'\mid s,a)
r_t(s,a,s')
\bar y_t(s,a)
+
\sum_{s\in\mathbb{S}}
r_T(s)\bar y_T(s)\\
\geq\;&
\big(
1-\delta
\big)
\sum_{t\in\mathbb{T}}
\sum_{s\in\mathbb{S}}
\sum_{a\in\mathcal{A}_t(s)}
\sum_{s'\in\mathbb{S}}
P_t(s'\mid s,a)
r_t(s,a,s')
\bar y_t(s,a)
+
\big(
1-\delta
\big)
\sum_{s\in\mathbb{S}}
r_T(s)\bar y_T(s)
=
\big(
1-\delta
\big)
F(\bar{\mathbf{y}}),
\end{aligned}
\end{equation*}
where the second inequality follows from $r_T(s)\geq0$,
$\bar y_T(s)\geq0$, and $\delta\in(0,1)$. Let
$\mathbf{y}^*\in\mathcal{Y}$ be an optimal solution to the exact
dual formulation \eqref{MDP_dual_LP}. Since $\mathbf{y}^*$ is also feasible for the dual
formulation \eqref{MDP_dual_LP} based on the local incumbent revenues, we have $\bar Z_{\mathrm{D}}
=
\max_{\mathbf{y}\in\mathcal{Y}}
\bar F(\mathbf{y})
\geq
\bar F(\mathbf{y}^*)
\geq
\big(
1-\delta
\big)
F(\mathbf{y}^*)
=
\big(
1-\delta
\big)
Z_{\mathrm{D}}$. Combining this relation with the baseline relation gives
$
\big(
1-\delta
\big)
Z_{\mathrm{D}}
\leq
\bar Z_{\mathrm{D}}
\leq
Z_{\mathrm{D}}.
$ Subtracting $\bar Z_{\mathrm{D}}$ from $Z_{\mathrm{D}}$ gives the
equivalent relation
$
0
\leq
Z_{\mathrm{D}}
-
\bar Z_{\mathrm{D}}
\leq
\delta Z_{\mathrm{D}}.
$
\end{proof}

Proposition
\ref{prop:dual_uniform_local_optimality_gap} uses the largest local optimality gap $\delta$ to analyze the sensitivity of the dual optimal value to changes in the local revenues when the transition and terminal revenues are nonnegative. We can also retain the incumbent revenue and upper bound obtained for each local problem and
use the transition probabilities to evaluate the effect of the individual local optimality gaps on the dual optimal value. Here, the interval of 
$[\bar r_t(s,a,s'),U_t(s,a,s')]$ describes the possible variation in a local revenue coefficient. We compute the largest expected change in a dual objective coefficient at each epoch and sum these changes over the horizon. That is, we define
\begin{equation*}
\begin{aligned}
\Delta_t
=
\max_{\substack{
s\in\mathbb{S}\\
a\in\mathcal{A}_t(s)}}
\left\{
\sum_{s'\in\mathbb{S}}
P_t(s'\mid s,a)
\big(
U_t(s,a,s')
-
\bar r_t(s,a,s')
\big)
\right\}
\quad 
\forall t\in\mathbb{T}, \text{ and }
\Delta_{\mathrm{D}}
=
\sum_{t\in\mathbb{T}}
\Delta_t.
\end{aligned}
\end{equation*}
Based on the definition of the local optimality gap, if $\delta_t(s,a,s')$ is bounded, then $U_t(s,a,s')
-
\bar r_t(s,a,s')
=
\delta_t(s,a,s')
\left|
U_t(s,a,s')
\right|$.
If $\bar Z_{\mathrm{D}}>0$, we define $\delta_{\mathrm{D}}
= \Delta_{\mathrm{D}}/ (\bar Z_{\mathrm{D}}
+
\Delta_{\mathrm{D}})$. The value $\delta_{\mathrm{D}}$ provides the total change
$\Delta_{\mathrm{D}}$ relative to
$\bar Z_{\mathrm{D}}+\Delta_{\mathrm{D}}$. The following proposition
provides an alternative sensitivity analysis. Proposition
\ref{prop:dual_uniform_local_optimality_gap} uses only the largest local optimality gap $\delta$. The result below also retains the individual local revenue intervals and transition probabilities through $\Delta_{\mathrm{D}}$ and compares $\delta_{\mathrm{D}}$ with $\delta$.

\begin{proposition}
\label{prop:dual_heterogeneous_local_optimality_bounds}
Suppose that $\bar Z_{\mathrm{D}}>0$, the transition and
terminal revenues are nonnegative, and there exists
$\delta\in(0,1)$ such that $0
\leq
\delta_t(s,a,s')
\leq
\delta$ for all $a\in\mathcal{A}_t(s),\
s,s'\in\mathbb{S},\
t\in\mathbb{T}.$ Then, we have
$0
\leq
Z_{\mathrm{D}}
-
\bar Z_{\mathrm{D}}
\leq
\min
\left\{
\delta,
\delta_{\mathrm{D}}
\right\}
Z_{\mathrm{D}}.$
\end{proposition}
\begin{proof}
Let $F(\mathbf{y})$ and $\bar F(\mathbf{y})$ denote the
objective functions of the exact dual formulation and the dual
formulation based on the local incumbent revenues, respectively.
Based on the local incumbent and upper bound relation, we have $0
\leq
r_t(s,a,s')
-
\bar r_t(s,a,s')
\leq
U_t(s,a,s')
-
\bar r_t(s,a,s')$
for all $t\in\mathbb{T}$, $s,s'\in\mathbb{S}$, and
$a\in\mathcal{A}_t(s)$. Therefore,
$\Delta_t\geq0$ for every $t\in\mathbb{T}$ and
$\Delta_{\mathrm{D}}\geq0$.

Let $\mathbf{y}^*\in\mathcal{Y}$ be an optimal solution to the
exact dual formulation. Based on the baseline relation
$\bar Z_{\mathrm{D}}\leq Z_{\mathrm{D}}$ and the feasibility of
$\mathbf{y}^*$ for the dual formulation based on the local
incumbent revenues, we have
\begin{equation*}
\begin{aligned}
0
\leq
Z_{\mathrm{D}}
-
\bar Z_{\mathrm{D}}
&=
F(\mathbf{y}^*)
-
\max_{\mathbf{y}\in\mathcal{Y}}
\bar F(\mathbf{y})\leq
F(\mathbf{y}^*)
-
\bar F(\mathbf{y}^*)\\
&=
\sum_{t\in\mathbb{T}}
\sum_{s\in\mathbb{S}}
\sum_{a\in\mathcal{A}_t(s)}
y_t^*(s,a)
\sum_{s'\in\mathbb{S}}
P_t(s'\mid s,a)
\left(
r_t(s,a,s')
-
\bar r_t(s,a,s')
\right)\\
&\leq
\sum_{t\in\mathbb{T}}
\sum_{s\in\mathbb{S}}
\sum_{a\in\mathcal{A}_t(s)}
y_t^*(s,a)
\sum_{s'\in\mathbb{S}}
P_t(s'\mid s,a)
\left(
U_t(s,a,s')
-
\bar r_t(s,a,s')
\right)\\
&\leq
\sum_{t\in\mathbb{T}}
\Delta_t
\underbrace{
\sum_{s\in\mathbb{S}}
\sum_{a\in\mathcal{A}_t(s)}
y_t^*(s,a)
}_{=\,1}=
\sum_{t\in\mathbb{T}}
\Delta_t
=
\Delta_{\mathrm{D}}.
\end{aligned}
\end{equation*}
The first upper bound follows because $\mathbf{y}^*$ is feasible
for the dual formulation based on the local incumbent revenues. The
first equality after this bound follows because both dual
formulations use the same terminal revenues. The next inequality
follows from
$r_t(s,a,s')\leq U_t(s,a,s')$ and
$P_t(s'\mid s,a)y_t^*(s,a)\geq0$. The last inequality follows from
the definition of $\Delta_t$ and Lemma
\ref{lem:dual_probability_flow_bounds}.

Since $\bar Z_{\mathrm{D}}>0$ and
$\Delta_{\mathrm{D}}\geq0$, the definition of
$\delta_{\mathrm{D}}$ gives $1
-
\delta_{\mathrm{D}}
= \bar Z_{\mathrm{D}}/(\bar Z_{\mathrm{D}}
+
\Delta_{\mathrm{D}}) > 0$. Since $Z_{\mathrm{D}}-\bar Z_{\mathrm{D}}\leq\Delta_{\mathrm{D}}$, adding $\bar Z_{\mathrm{D}}$ to both sides gives $Z_{\mathrm{D}}\leq\bar Z_{\mathrm{D}}+\Delta_{\mathrm{D}}$. Since
$Z_{\mathrm{D}}
\leq
\bar Z_{\mathrm{D}}+\Delta_{\mathrm{D}}$
and
$\bar Z_{\mathrm{D}}/
(\bar Z_{\mathrm{D}}+\Delta_{\mathrm{D}})>0$, multiplying both sides
by
$\bar Z_{\mathrm{D}}/
(\bar Z_{\mathrm{D}}+\Delta_{\mathrm{D}})$
gives $ (\bar Z_{\mathrm{D}} / (\bar Z_{\mathrm{D}} + \Delta_{\mathrm{D}})) Z_{\mathrm{D}} \leq (\bar Z_{\mathrm{D}} / (\bar Z_{\mathrm{D}} + \Delta_{\mathrm{D}})) (\bar Z_{\mathrm{D}} + \Delta_{\mathrm{D}}) = \bar Z_{\mathrm{D}} $. Based on
$
1-\delta_{\mathrm{D}}
=
\bar Z_{\mathrm{D}}/
(\bar Z_{\mathrm{D}}+\Delta_{\mathrm{D}})
$,
this relation becomes $\big(
1-\delta_{\mathrm{D}}
\big)
Z_{\mathrm{D}}
\leq
\bar Z_{\mathrm{D}}$. The inequality
$\big(1-\delta_{\mathrm{D}}\big)Z_{\mathrm{D}}
\leq\bar Z_{\mathrm{D}}$
gives
$Z_{\mathrm{D}}-\bar Z_{\mathrm{D}}
\leq\delta_{\mathrm{D}}Z_{\mathrm{D}}$.
The baseline relation
$\bar Z_{\mathrm{D}}\leq Z_{\mathrm{D}}$
gives
$0\leq Z_{\mathrm{D}}-\bar Z_{\mathrm{D}}$.
Combining these two inequalities gives
$0\leq Z_{\mathrm{D}}-\bar Z_{\mathrm{D}}
\leq\delta_{\mathrm{D}}Z_{\mathrm{D}}$. Based on Proposition
\ref{prop:dual_uniform_local_optimality_gap}, we also have
$0\leq Z_{\mathrm{D}}-\bar Z_{\mathrm{D}}
\leq\delta Z_{\mathrm{D}}$. Finally, combining these two bounds gives $0
\leq
Z_{\mathrm{D}}
-
\bar Z_{\mathrm{D}}
\leq
\min
\left\{
\delta,
\delta_{\mathrm{D}}
\right\}
Z_{\mathrm{D}}.$
This completes the proof.
\end{proof}

Propositions
\ref{prop:dual_uniform_local_optimality_gap} and
\ref{prop:dual_heterogeneous_local_optimality_bounds} analyze the sensitivity of the MDP optimal value to local optimality gaps under nonnegative transition and terminal revenues. Proposition
\ref{prop:dual_heterogeneous_local_optimality_bounds} additionally assumes that $\bar Z_{\mathrm{D}}>0$. The following two examples satisfy these conditions. They provide two sensitivity evaluations for the MDP optimal value. One example uses
the largest local optimality gap, and the other uses the individual
incumbent revenues and upper bounds. The examples also compare the
resulting evaluations of the global MDP optimality gap.

\begin{example}
Suppose that we have an instance with $T=2$, $\mathbb{S}=\{s\}$,
$\mathcal{A}_t(s)=\{a\}$, $P_t(s\mid s,a)=1$ for
$t\in\{0,1\}$, and $r_T(s)=0$. Besides, the local incumbent
revenues, upper bounds, and relative optimality gaps are
\begin{equation*}
\begin{aligned}
\big(
\bar r_0(s,a,s),
U_0(s,a,s),
\delta_0(s,a,s)
\big)
=
(99,100,0.01) \text{ and }
\big(
\bar r_1(s,a,s),
U_1(s,a,s),
\delta_1(s,a,s)
\big)
=
(8,10,0.20).
\end{aligned}
\end{equation*}
The intervals $[99,100]$ and $[8,10]$ describe the possible changes
in the local revenue coefficients. Then, we have
$\bar Z_{\mathrm{D}}
=
99+8
=
107$ and 
$\delta
=
\max
\left\{
0.01,
0.20
\right\}
=
0.20.$
Based on the relation
$0\leq
Z_{\mathrm{D}}-\bar Z_{\mathrm{D}}
\leq
\delta Z_{\mathrm{D}}$
of Proposition
\ref{prop:dual_uniform_local_optimality_gap}, we have
$0
\leq
Z_{\mathrm{D}}
-
107
\leq
0.20Z_{\mathrm{D}}$.
Thus, this sensitivity analysis gives the relations
$0.8Z_{\mathrm{D}}
\leq
107
\leq
Z_{\mathrm{D}}$ and $107
\leq
Z_{\mathrm{D}}
\leq
133.75$.
For the sensitivity analysis of Proposition
\ref{prop:dual_heterogeneous_local_optimality_bounds}, we have
\begin{equation*}
\Delta_0
=
100-99
=
1,
\Delta_1
=
10-8
=
2,
\Delta_{\mathrm{D}}
=
3,
\text{and }
\delta_{\mathrm{D}}
=
\frac{
\Delta_{\mathrm{D}}
}{
\bar Z_{\mathrm{D}}
+
\Delta_{\mathrm{D}}
}
=
\frac{3}{107+3}
=
\frac{3}{110}
\approx
0.0273.
\end{equation*}
Based on the relation
$0\leq
Z_{\mathrm{D}}-\bar Z_{\mathrm{D}}
\leq
\min\{\delta,\delta_{\mathrm{D}}\}Z_{\mathrm{D}}$
of Proposition
\ref{prop:dual_heterogeneous_local_optimality_bounds}, we have $0
\leq
Z_{\mathrm{D}}
-
107
\leq
\min\{0.2,\frac{3}{110}\}Z_{\mathrm{D}}=
\frac{3}{110}
Z_{\mathrm{D}}$ and $107
\leq
Z_{\mathrm{D}}
\leq
110.$
Thus, the sensitivity analysis using the incumbent revenue and upper
bound of each local problem in Proposition
\ref{prop:dual_heterogeneous_local_optimality_bounds} limits the
possible relative change in the dual optimal value to approximately
$2.73\%$, compared with $20\%$ obtained from the largest local
optimality gap. If
$r_0(s,a,s)=100$ and $r_1(s,a,s)=10$, then
$Z_{\mathrm{D}}=110$, and the latter bound is attained.
\end{example}

We note that Proposition
\ref{prop:dual_uniform_local_optimality_gap} has a broader range of
applicability because its sensitivity result remains valid for
$\bar Z_{\mathrm{D}}\geq0$, including the case
$\bar Z_{\mathrm{D}}=0$. Proposition
\ref{prop:dual_heterogeneous_local_optimality_bounds} requires
$\bar Z_{\mathrm{D}}>0$ to define $\delta_{\mathrm{D}}$. When both
propositions apply, we have $\min
\left\{
\delta,
\delta_{\mathrm{D}}
\right\}
\leq
\delta.$ Therefore, the combined sensitivity measure in Proposition
\ref{prop:dual_heterogeneous_local_optimality_bounds} is no larger
than the measure in Proposition
\ref{prop:dual_uniform_local_optimality_gap}. The ordering of the
individual sensitivity measures $\delta$ and
$\delta_{\mathrm{D}}$ depends on $\Delta_{\mathrm{D}}$ and
$\bar Z_{\mathrm{D}}$. Based on the definition of
$\delta_{\mathrm{D}}$, the following relations
\begin{equation*}
\delta_{\mathrm{D}}
\leq
\delta,
\quad
\frac{
\Delta_{\mathrm{D}}
}{
\bar Z_{\mathrm{D}}
+
\Delta_{\mathrm{D}}
}
\leq
\delta,
\quad
\big(
1-\delta
\big)
\Delta_{\mathrm{D}}
\leq
\delta\bar Z_{\mathrm{D}},
\quad
\text{and}
\quad
\Delta_{\mathrm{D}}
\leq
\frac{
\delta
}{
1-\delta
}
\bar Z_{\mathrm{D}}
\end{equation*}
are equivalent.
Thus, we can observe that $\delta_{\mathrm{D}}<\delta$ if
$\Delta_{\mathrm{D}}
<
\delta\bar Z_{\mathrm{D}}/(1-\delta)$, and
$\delta<\delta_{\mathrm{D}}$ if
$\Delta_{\mathrm{D}}
>
\delta\bar Z_{\mathrm{D}}/(1-\delta)$. Therefore, neither
$\delta$ nor $\delta_{\mathrm{D}}$ gives the smaller individual
sensitivity measure for every instance. The following example shows
that the common upper bound $\delta$ on the local relative optimality
gaps can instead give the smaller sensitivity measure.

\begin{example}
Suppose that we have an instance with $T=2$. At $t=0$, the initial state $s_0$
has one action that reaches $s_1$ and $s_2$ with probabilities
$0.999$ and $0.001$, respectively. Let all transition revenues,
local incumbent revenues, and upper bounds at $t=0$ be zero, and let
$r_T(s)=0$ for all $s\in\mathbb{S}$. At $t=1$, each state has one
action and a deterministic transition, with
\begin{equation*}
\begin{aligned}
\big(
r_1(s_1,a_1,s_1),
\bar r_1(s_1,a_1,s_1),
U_1(s_1,a_1,s_1)
\big)
&=
(100,90,100),\text{ and}\\
\big(
r_1(s_2,a_2,s_2),
\bar r_1(s_2,a_2,s_2),
U_1(s_2,a_2,s_2)
\big)
&=
(100000,90000,100000).
\end{aligned}
\end{equation*}
The intervals $[90,100]$ and $[90000,100000]$ describe the possible
changes in the two local revenue coefficients. Then, we have
$\delta_1(s_1,a_1,s_1)
=
\delta_1(s_2,a_2,s_2)
=
\delta
=
0.10$,
$Z_{\mathrm{D}}
=
0.999(100)+0.001(100000)
=
199.9$, and
$\bar Z_{\mathrm{D}}
=
0.999(90)+0.001(90000)
=
179.91$.
Based on the relation
$0\leq
Z_{\mathrm{D}}-\bar Z_{\mathrm{D}}
\leq
\delta Z_{\mathrm{D}}$
of Proposition
\ref{prop:dual_uniform_local_optimality_gap}, we have
$0\leq
Z_{\mathrm{D}}-\bar Z_{\mathrm{D}}
\leq
0.10Z_{\mathrm{D}}$. For the sensitivity analysis of Proposition
\ref{prop:dual_heterogeneous_local_optimality_bounds}, we have
\begin{equation*}
\begin{aligned}
\Delta_0
=
0,
\Delta_1
=
\max
\left\{
100-90,\,
100000-90000
\right\}
=
10000,
\Delta_{\mathrm{D}}
=
10000,
\text{ and }
\delta_{\mathrm{D}}
=
\frac{
10000
}{
179.91+10000
}
\approx
0.9823.
\end{aligned}
\end{equation*}
Here, $\Delta_1$ is determined by the larger local revenue interval
at $s_2$, although $s_2$ is reached with probability $0.001$. Based
on the relation
$0\leq
Z_{\mathrm{D}}-\bar Z_{\mathrm{D}}
\leq
\min\{\delta,\delta_{\mathrm{D}}\}Z_{\mathrm{D}}$
of Proposition
\ref{prop:dual_heterogeneous_local_optimality_bounds}, we have
$0\leq
Z_{\mathrm{D}}-\bar Z_{\mathrm{D}}
\leq
\min\{0.10,0.9823\}Z_{\mathrm{D}}
=
0.10Z_{\mathrm{D}}$.
Thus, the common upper bound on the local relative optimality gaps in
Proposition
\ref{prop:dual_uniform_local_optimality_gap}
gives the smaller individual sensitivity measure in this instance.
\end{example}

The above examples compare the sensitivity measures obtained
from the largest local optimality gap and the individual local
revenue intervals. We next use the second sensitivity measure to
develop a stopping condition for the local solution process. An
increase in a local incumbent revenue cannot decrease
$\bar Z_{\mathrm{D}}$, and a decrease in a local upper bound cannot
increase $\Delta_{\mathrm{D}}$. Therefore, the local incumbent
revenues and upper bounds can be updated until a target global
relative optimality gap is satisfied. In particular, the target
global optimality gap may be satisfied even if some local optimality
gaps remain larger than the target value. We provide the following corollary that states this stopping condition.
\begin{corollary}
\label{cor:dual_global_stopping_condition}
Under the conditions of Proposition
\ref{prop:dual_heterogeneous_local_optimality_bounds}, let
$\eta\in(0,1)$ be a target global relative optimality gap. If $\delta
\leq
\eta
\text{ or }
\Delta_{\mathrm{D}}
\leq
\frac{
\eta
}{
1-\eta
}
\bar Z_{\mathrm{D}}$, then $0
\leq
Z_{\mathrm{D}}
-
\bar Z_{\mathrm{D}}
\leq
\eta Z_{\mathrm{D}}.$
In particular, the second condition does not require
$\delta\leq\eta$.
\end{corollary}
\begin{proof}
If $\delta\leq\eta$, the result follows directly from Proposition
\ref{prop:dual_heterogeneous_local_optimality_bounds}. If $\Delta_{\mathrm{D}}
\leq
\frac{
\eta
}{
1-\eta
}
\bar Z_{\mathrm{D}},$ then $\big(
1-\eta
\big)
\Delta_{\mathrm{D}}
\leq
\eta\bar Z_{\mathrm{D}}$ and $\Delta_{\mathrm{D}}
\leq
\eta
\left(
\bar Z_{\mathrm{D}}
+
\Delta_{\mathrm{D}}
\right)$.
This result gives $\delta_{\mathrm{D}}
= \Delta_{\mathrm{D}}/(\bar Z_{\mathrm{D}}
+
\Delta_{\mathrm{D}})
\leq
\eta.$
Thus,
$\min\{\delta,\delta_{\mathrm{D}}\}\leq\eta$, and Proposition
\ref{prop:dual_heterogeneous_local_optimality_bounds} gives
$0
\leq
Z_{\mathrm{D}}
-
\bar Z_{\mathrm{D}}
\leq
\eta Z_{\mathrm{D}}.$
\end{proof}
The following example considers two actions at the initial state. One
action leads to a successor state with a local revenue problem having
a positive optimality gap. The other action leads to a successor
state with a known revenue. We continue solving only the first local
problem and update its incumbent revenue. The example shows that the
second condition of Corollary
\ref{cor:dual_global_stopping_condition} can satisfy a target global
relative optimality gap even when the remaining local optimality gap
is larger than the target value.
\begin{example}
Suppose that we have an instance with $T=2$,
$\mathbb{S}=\{s_0,s^{\mathrm{I}},s^{\mathrm{N}}\}$,
$\mathcal{A}_0(s_0)=\{a^{\mathrm{I}},a^{\mathrm{N}}\}$, and
$\mathcal{A}_1(s^{\mathrm{I}})
=
\mathcal{A}_1(s^{\mathrm{N}})
=
\{a\}$, where the superscripts $\mathrm{I}$ and $\mathrm{N}$
represent installation and no installation, respectively. Besides, we have 
$P_0(s^{\mathrm{I}}\mid s_0,a^{\mathrm{I}})
=
P_0(s^{\mathrm{N}}\mid s_0,a^{\mathrm{N}})
=
P_1(s^{\mathrm{I}}\mid s^{\mathrm{I}},a)
=
P_1(s^{\mathrm{N}}\mid s^{\mathrm{N}},a)
=
1$ and $r_T(s)=0$ for all $s\in\mathbb{S}$. The local incumbent
revenues, exact revenues, upper bounds, and relative optimality gaps
are
\begin{equation*}
\begin{aligned}
\big(
\bar r_0(s_0,a^{\mathrm{I}},s^{\mathrm{I}}),
r_0(s_0,a^{\mathrm{I}},s^{\mathrm{I}}),
U_0(s_0,a^{\mathrm{I}},s^{\mathrm{I}})
\big)
&=
(60,60,60),\text{ and}\\
\big(
\bar r_0(s_0,a^{\mathrm{N}},s^{\mathrm{N}}),
r_0(s_0,a^{\mathrm{N}},s^{\mathrm{N}}),
U_0(s_0,a^{\mathrm{N}},s^{\mathrm{N}})
\big)
&=
(100,100,100),\\
\big(
\bar r_1(s^{\mathrm{I}},a,s^{\mathrm{I}}),
U_1(s^{\mathrm{I}},a,s^{\mathrm{I}}),
\delta_1(s^{\mathrm{I}},a,s^{\mathrm{I}})
\big)
&=
(50,100,0.50),\text{ and}\\
\big(
\bar r_1(s^{\mathrm{N}},a,s^{\mathrm{N}}),
r_1(s^{\mathrm{N}},a,s^{\mathrm{N}}),
U_1(s^{\mathrm{N}},a,s^{\mathrm{N}})
\big)
&=
(0,0,0).
\end{aligned}
\end{equation*}
Then, we have
$\bar Z_{\mathrm{D}}
=
\max\{60+50,100\}
=
110$,
$\delta=0.50$,
$\Delta_0=0$,
$\Delta_1
=
\max\{100-50,0\}
=
50$,
$\Delta_{\mathrm{D}}=50$, and
$\delta_{\mathrm{D}}
=
50/(110+50)
=
0.3125$.
For the target global relative optimality gap $\eta=0.25$, neither
condition of Corollary
\ref{cor:dual_global_stopping_condition} holds because
\begin{equation*}
\delta
=
0.50
>
0.25
\quad\text{and}\quad
\Delta_{\mathrm{D}}
=
50
>
\frac{0.25}{1-0.25}
110
=
\frac{110}{3}.
\end{equation*}
Suppose that the local problem associated with
$(1,s^{\mathrm{I}},a,s^{\mathrm{I}})$ is solved further and its
incumbent revenue increases from $50$ to $60$, while its upper bound
remains equal to $100$. Then, we have
\begin{equation*}
\big(
\bar r_1^{+}(s^{\mathrm{I}},a,s^{\mathrm{I}}),
U_1^{+}(s^{\mathrm{I}},a,s^{\mathrm{I}}),
\delta_1^{+}(s^{\mathrm{I}},a,s^{\mathrm{I}})
\big)
=
(60,100,0.40).
\end{equation*}
The updated values satisfy
$\bar Z_{\mathrm{D}}^{+}
=
\max\{60+60,100\}
=
120$,
$\delta^{+}=0.40$,
$\Delta_{\mathrm{D}}^{+}
=
100-60
=
40$, and
$\delta_{\mathrm{D}}^{+}
=
40/(120+40)
=
0.25$.
Although $\delta^{+}=0.40>0.25$, we have
\begin{equation*}
\Delta_{\mathrm{D}}^{+}
=
40
=
\frac{0.25}{1-0.25}
120
=
\frac{\eta}{1-\eta}
\bar Z_{\mathrm{D}}^{+}.
\end{equation*}
Based on the second condition of Corollary
\ref{cor:dual_global_stopping_condition}, we have
$0
\leq
Z_{\mathrm{D}}
-
\bar Z_{\mathrm{D}}^{+}
\leq
0.25Z_{\mathrm{D}}$.
Thus, the local solution process can terminate with a $40\%$ local
optimality gap even though the target global relative optimality gap
is $25\%$.
\end{example}

The above sensitivity results from Propositions
\ref{prop:dual_uniform_local_optimality_gap} and
\ref{prop:dual_heterogeneous_local_optimality_bounds} use $\delta$ and
$\Delta_{\mathrm{D}}$ to evaluate changes in the local revenue
coefficients. Finally, we set every local revenue coefficient to the
upper endpoint of its interval and retain the original dual feasible
region $\mathcal{Y}$. We define
\begin{equation*}
Z_{\mathrm{D}}^{U}
=
\max_{\mathbf{y}\in\mathcal{Y}}
\left\{
\sum_{t\in\mathbb{T}}
\sum_{s\in\mathbb{S}}
\sum_{a\in\mathcal{A}_t(s)}
\sum_{s'\in\mathbb{S}}
P_t(s'\mid s,a)
U_t(s,a,s')
y_t(s,a)
+
\sum_{s\in\mathbb{S}}
r_T(s)y_T(s)
\right\}.
\end{equation*}
Since the local revenue coefficients do not affect the dual feasible
region, $Z_{\mathrm{D}}^{U}$ is the largest dual optimal value allowed
by the local revenue intervals
$[\bar r_t(s,a,s'),U_t(s,a,s')]$. If
$\bar Z_{\mathrm{D}}>0$, we define the corresponding relative
sensitivity gap as $\delta_{\mathrm{D}}^{U} = (Z_{\mathrm{D}}^{U}
-
\bar Z_{\mathrm{D}})/Z_{\mathrm{D}}^{U}$.
The following corollary concludes this sensitivity analysis. It
relates $Z_{\mathrm{D}}^{U}-\bar Z_{\mathrm{D}}$ to
$\Delta_{\mathrm{D}}$ and compares $\delta$,
$\delta_{\mathrm{D}}$, and $\delta_{\mathrm{D}}^{U}$.
\begin{corollary}
\label{cor:dual_reoptimized_local_upper_bounds}
Given
$\bar r_t(s,a,s')
\leq
r_t(s,a,s')
\leq
U_t(s,a,s')$
for all $t\in\mathbb{T}$, $s,s'\in\mathbb{S}$, and
$a\in\mathcal{A}_t(s)$, we have
$0
\leq
Z_{\mathrm{D}}
-
\bar Z_{\mathrm{D}}
\leq
Z_{\mathrm{D}}^{U}
-
\bar Z_{\mathrm{D}}
\leq
\Delta_{\mathrm{D}}$.
If the conditions of Proposition
\ref{prop:dual_heterogeneous_local_optimality_bounds} hold, then
$0
\leq
\frac{
Z_{\mathrm{D}}
-
\bar Z_{\mathrm{D}}
}{
Z_{\mathrm{D}}
}
\leq
\delta_{\mathrm{D}}^{U}
\leq
\min
\left\{
\delta,
\delta_{\mathrm{D}}
\right\}
<
1$.
\end{corollary}
\begin{proof}
For every $\mathbf{y}\in\mathcal{Y}$, the local incumbent and upper
bound relations give
$\bar F(\mathbf{y})
\leq
F(\mathbf{y})
\leq
F^{U}(\mathbf{y})$.
Let $\mathbf{y}^{U}\in\mathcal{Y}$ be an optimal solution associated
with $Z_{\mathrm{D}}^{U}$. Based on the definition of
$\Delta_{\mathrm{D}}$ and Lemma
\ref{lem:dual_probability_flow_bounds}, we have $0
\leq
Z_{\mathrm{D}}
-
\bar Z_{\mathrm{D}}
\leq
Z_{\mathrm{D}}^{U}
-
\bar Z_{\mathrm{D}}
\leq
F^{U}(\mathbf{y}^{U})
-
\bar F(\mathbf{y}^{U})
\leq
\Delta_{\mathrm{D}}.$
Under the conditions of Proposition
\ref{prop:dual_heterogeneous_local_optimality_bounds}, we have
$\bar Z_{\mathrm{D}}>0$ and $\frac{
Z_{\mathrm{D}}
-
\bar Z_{\mathrm{D}}
}{
Z_{\mathrm{D}}
}
\leq
\frac{
Z_{\mathrm{D}}^{U}
-
\bar Z_{\mathrm{D}}
}{
Z_{\mathrm{D}}^{U}
}
=
\delta_{\mathrm{D}}^{U}
\leq
\delta_{\mathrm{D}}.$
Moreover,
$\bar r_t(s,a,s')
\geq
(1-\delta)U_t(s,a,s')$
and the nonnegative terminal revenues give
$\bar Z_{\mathrm{D}}
\geq
(1-\delta)Z_{\mathrm{D}}^{U}$.
Thus, $\delta_{\mathrm{D}}^{U}\leq\delta$. This 
completes the proof.
\end{proof}
We conclude this section with a summary of the sensitivity
analyses. Proposition
\ref{prop:dual_uniform_local_optimality_gap} uses the largest local
optimality gap $\delta$ to evaluate a common relative change in the
local revenue coefficients. Proposition
\ref{prop:dual_heterogeneous_local_optimality_bounds} uses the
individual local revenue intervals and transition probabilities to
derive $\Delta_{\mathrm{D}}$ and $\delta_{\mathrm{D}}$. Corollary
\ref{cor:dual_reoptimized_local_upper_bounds} uses all local upper
bounds over the original dual feasible region $\mathcal{Y}$ to
derive $Z_{\mathrm{D}}^{U}$ and
$\delta_{\mathrm{D}}^{U}$. These results provide three levels of
sensitivity analysis based on $\delta$,
$\delta_{\mathrm{D}}$, and
$\delta_{\mathrm{D}}^{U}$, respectively. 

\section{Conclusion}

Motivated by the stochastic depletion framework of Chan and Farias
\cite{ChanMOR2009}, we consider a finite horizon MDP model \eqref{MDP_model} with transition
revenues defined by the local submodular revenue maximization problems
\eqref{MDP_transition_revenue}. We study an initial mixed-integer
programming analysis of this structure. We observe that
the optimal local revenues appear as coefficients in the objective
function of the dual formulation \eqref{MDP_dual_LP}, and the feasible
region of the MDP dual variables remains unchanged. Therefore, the
dual formulation provides a connection between the MDP variables and
the binary variables of the local problems.

Based on this dual structure, we obtain an exact mixed-integer
nonlinear formulation with bilinear terms. We then use the
submodular inequalities and an exact McCormick linearization to
derive an equivalent mixed-integer linear formulation.  The resulting MILP provides a basis for further polyhedral analysis
of the combined structure defined by
$\mathcal{D}^{(t,s,a,s')}$ and
$\mathcal{W}^{(t,s,a,s')}$. In addition, the MDP model \eqref{MDP_model} contains many NP-hard local
submodular revenue maximization problems
\eqref{MDP_transition_revenue}. We show that the dual formulation provides a basis
for analyzing the sensitivity of the MDP optimal value to changes in
the local revenues. Our sensitivity analysis uses the incumbent revenues and valid
upper bounds of the local problems to relate the local optimality gaps
to the global MDP optimality gap.

\section*{Acknowledgments}
This work is supported, in part by, NSTC Taiwan 114-2221-E-A49-176-MY2. All authors are listed in alphabetical order.

	\bibliographystyle{abbrv}
	\bibliography{general}

\end{document}